\documentclass[b4paper,9pt]{article}
\usepackage{array}

\RequirePackage[OT1]{fontenc} \RequirePackage{amsmath,natbib}
\RequirePackage[colorlinks,citecolor=black,urlcolor=blue]{hyperref}
\RequirePackage{hypernat}
\usepackage{enumerate}
\usepackage{amsfonts}
\usepackage{graphicx}
\usepackage{verbatim}
\usepackage{amsthm}
\usepackage{url}
\usepackage{bm}
\usepackage{natbib}
\usepackage{amsmath}
\usepackage{extarrows}
\usepackage{amsfonts}
\usepackage{amssymb}

\newtheorem{definition}{Definition}
\newtheorem{theorem}{Theorem}
\newtheorem{corollary}{Corollary}
\newtheorem{lemma}{Lemma}
\newtheorem{proposition}{Proposition}
\newtheorem{remark}{Remark}
\newtheorem{example}{Example}

\newcommand{\bfw}{\mathbf{w}}

\newcommand{\E}{\mathrm{E}}
\newcommand{\pr}{\mathrm{P}}
\newcommand{\R}{\mathbb{R}}

\newcommand{\hl}[1]{{#1}}

\usepackage{geometry}
\begin{document}


\renewcommand{\baselinestretch}{1}

\markright{
\hbox{\footnotesize\rm Preprint}\hfill
}

\markboth{\hfill{\footnotesize\rm ~~~} \hfill}
{\hfill {\footnotesize\rm ~~~} \hfill}

\renewcommand{\thefootnote}{}
$\ $\par


\fontsize{10.95}{14pt plus.8pt minus .6pt}\selectfont
\vspace{0.8pc}
\centerline{\large\bf Sharper Sub-Weibull Concentrations}
\vspace{.4cm}
\centerline{Huiming Zhang $^{1,2}$ and Haoyu Wei$^{3 \dag}$\footnote{* Two authors contribute equally to this work. Huiming Zhang's e-mail: huimingzhang@um.edu.mo; Haoyu Wei's e-mail: cute@pku.edu.cn} \footnote{\dag \,  Corresponding author.}}
\vspace{.4cm}
\centerline{\it 1. Department of Mathematics,
University of Macau, Taipa Macau, China.}
~\\
\centerline{\it 2.  UMacau Zhuhai Research Institute, Zhuhai, China.}
~\\
\centerline{\it 3. Guanghua School of Management, Peking University, Beijing, China.}
\vspace{.55cm}
\fontsize{9}{11.5pt plus.8pt minus .6pt}\selectfont

\begin{quotation}
\noindent {\it Abstract:}
Constant-specified and exponential concentration inequalities play an essential role in the finite-sample theory of machine learning and high-dimensional statistics area. We obtain sharper and constants-specified concentration inequalities for the sum of independent sub-Weibull random variables, which leads to a mixture of two tails: sub-Gaussian for small deviations and sub-Weibull for large deviations from the mean. These bounds are new and improve existing bounds with sharper constants. In addition, {a new}
 \emph{sub-Weibull parameter}
 is also proposed, which enables recovering the tight concentration inequality for a random variable (vector). For statistical applications, we give an $\ell_2$-error of estimated coefficients in negative binomial regressions when the heavy-tailed covariates are sub-Weibull distributed with sparse structures, which is a new result for negative binomial regressions. In applying random matrices, we derive non-asymptotic versions of Bai-Yin's theorem for sub-Weibull entries with exponential tail bounds. Finally, by demonstrating a sub-Weibull confidence region for a log-truncated Z-estimator without the second-moment condition, we discuss and define the \emph{sub-Weibull type robust estimator} for independent observations $\{X_i\}_{i=1}^{n}$ without exponential-moment conditions.

\par

\vspace{9pt}
\par
\end{quotation}\par


\fontsize{10.95}{14pt plus.8pt minus .6pt}\selectfont

\section{Introduction}
In the last two decades, with the development of modern data collection methods in science and techniques, scientists and engineers can access and load a huge number of variables in their experiments. Over hundreds of years, probability theory lays the mathematical foundation of statistics. Arising from data-driving problems, various recent statistics research advances also contribute to new and challenging probability problems for further study. For example, in recent years, the rapid development of high-dimensional statistics and machine learning have promoted the development of the probability theory and even pure mathematics, such as random matrices, large deviation inequalities, and geometric functional analysis, etc.; see  \cite{Vershynin18}. More importantly, the concentration inequality (CI) quantifies the concentration of measures that are at the heart of statistical machine learning. Usually, CI quantifies how a random variable (r.v.) $X$ deviates around its mean $\mathrm{E}X=:\mu$ by presenting as one-side or two-sided bounds for the tail probability of $X-\mu$
$$
\pr (X-\mu >t)~\text{or}~\pr(|X-\mu|>t) \leq \text{some small}~\delta,~\forall~t \ge 0.
$$

The classical statistical models are faced with fixed-dimensional variables only. However, contemporary data science motivates statisticians to pay more attention to studying ${p \times p}$ random Hessian matrices (or sample covariance matrices, \cite {bai2010}) with $p \to \infty$, arising from the likelihood functions of high-dimensional regressions with covariates in $\R^p$. When the model dimension increases with sample size, obtaining asymptotic results for the estimator is potentially more challenging than the fixed dimensional case. In statistical machine learning, concentration inequalities (large derivation inequalities) are essential in deriving non-asymptotic error bounds for the proposed estimator; see \cite{Wainwright19,Zhang20}. Over recent decades, researchers have developed remarkable results of matrix concentration inequalities, which focuses on non-asymptotic upper and lower bounds for the largest eigenvalue of a finite sum of random matrices. For a more fascinated introduction, please refer to the book \cite{Tropp2015}.

Motivated from sample covariance matrices, a random matrix is a specific matrix $\mathbf{A}_{p \times p}$ with its entries $A_{jk}$ drawn from some distributions. As $p \to \infty$, random matrix theory mainly focuses on studying the properties of the $p$ eigenvalues of $\mathbf{A}_{p \times p}$, which turn out to have some limit law. Several famous limit laws in random matrix theory are different from the CLT for the summation of independent random variables since the $p$ eigenvalues are dependent and interact with each other. For convergence in distribution, some pioneering works are the Wigner's semicircle law for some symmetric Gaussian matrices' eigenvalues, the Marchenko-Pastur law for Wishart distributed random matrices (sample covariance matrices), and the Tracy-Widom laws for the limit distribution for maximum eigenvalues in Wishart matrices. All these three laws can be regarded as the CLT of random matrix versions. Moreover, the limit law for the empirical spectral density is some circle distribution, which sheds light on the non-communicative behaviors of the random matrix, while the classic limit law in CLT is for normal distribution or infinite divisible distribution. For strong convergence, Bai-Yin's law complements the Marchenko-Pastur law, which asserts that almost surely convergence of the smallest and largest eigenvalue for a sample covariance matrix. The monograph \cite{bai2010} thoroughly introduces the limit law in random matrices.

This work aims to extend non-asymptotic results from sub-Gaussian to sub-Weibull in terms of exponential concentration inequalities with applications in count data regressions, random matrices, and robust estimators. The contributions are: 

\begin{itemize}
    \item[(i)] We review and present some new results for sub-Weibull r.v.s, including sharp concentration inequalities for weighted summations of independent sub-Weibull r.v.s and negative binomial r.v.s, which are useful in many statistical applications.
    
    \item[(ii)] Based on the generalized Bernstein-Orlicz norm, a sharper concentration for sub-Weibull summations is obtained in Theorem \ref{thm:SumNewOrliczVex}. Here we circumvent Stirling's approximation and derive the inequalities more subtly. As a result, the confidence interval based on our result is sharper and more accurate than that in \cite{Kuchibhotla18} (For example, see Remark 2) and \cite{Hao19} (see Proposition 1 with unknown constants) gave.
    
    \item[(iii)] By sharper sub-Weibull concentrations, we give two applications. First, from the proposed negative binomial concentration inequalities, we obtain the $O_P(\sqrt{{p}/{n}})$ (up to some log factors) estimation error for the estimated coefficients in negative binomial regressions under the increasing-dimensional framework $p=p_n$ and heavy-tailed covariates. Second, we provide a non-asymptotic Bai-Yin's theorem for sub-Weibull random matrices with exponential-decay high probability. 
    \item[(iv)] We propose a new \emph{sub-Weibull parameters}, which is enabled of recovering the tight concentration inequality for a single non-zero mean random vector. The simulation studies for estimating sub-Gaussian and sub-exponential parameters show these parameters could be estimated well.
    
   \item[(v)] We establish a unified non-asymptotic confidence region and the convergence rate for general log-truncated Z-estimator in Theorem \ref{eq:cantoniG}. Moreover, we define a sub-Weibull type estimator for a sequence of independent observations $\{X_i\}_{i=1}^{n}$ without the second-moment condition, beyond the definition of the sub-Gaussian estimator.
\end{itemize}

\section{Sharper Concentrations for Sub-Weibull Summation}

Concentration inequalities are powerful in high-dimensional statistical inference, and it can derive explicit non-asymptotic error bounds as a function of sample size, sparsity level, and dimension \citep{Wainwright19}. In this section, we present preparation results of concentration inequalities for sub-Weibull random variables.

\subsection{Properties of Sub-Weibull norm and Orlicz-type  norm}
In empirical process theory, sub-Weibull norm (or other Orlicz-type norms) is crucial to derive the tail probability for both single sub-Weibull random variable and summation of random variables (by using the Chernoff's inequality). A benefit of Orlicz-type norms is that the concentration does not need the zero mean assumption.
\begin{definition}[Sub-Weibull norm]
    For $\theta>0$, the sub-Weibull norm of $X$ is defined as
    \begin{eqnarray*}
        \|X\|_{\psi_{\theta}}:=\inf\{C\in(0, \infty): ~ \mathrm E[\exp(|X|^{\theta}/C^{\theta})]\leq 2\}.
    \end{eqnarray*}
\end{definition}

The $\|\cdot\|_{\psi_{\theta}}$ is also called the $\psi_{\theta}$-norm. We define $X$ as a \emph{sub-Weibull random variable} with \emph{index} $\theta$ if it has a bounded $\psi_{\theta}$-norm (denoted as $X \sim \operatorname{subW} (\theta)$). Actually, the sub-Weibull norm is a special case of Orlicz norms below.

\begin{definition}[Orlicz Norms]\label{def:IncOrliczNorm}
\emph{Let $g:[0, \infty) \to [0, \infty)$ be a non-decreasing convex function with $g(0) = 1$. The ``$g$-Orlicz norm'' of a real-valued r.v. $X$ is given by}
\begin{equation}\label{eq:AlphaOrliczNorm}
\|X\|_{g}:=\inf \{\eta>0: \E[g(|X| / \eta)] \leq 2\}.
\end{equation}
\end{definition}
Using exponential Markov's inequality, we have
\begin{equation}\label{eq:tail}
    \pr(|X| \geq t)={\pr}({g(|X |/ \|X\|_g)} \geq {g(t / \|X\|_g)}) \leq g^{-1}(t / \|X\|_g) \mathrm{E} {g(X / \|X\|_g)} \leq 2 g^{-1}(t / \|X\|_g)
\end{equation}
by Definition \ref{def:IncOrliczNorm}. For example, let $g(x)=e^{x^{\theta}}$, which leads to sub-Weibull norm for $\theta \ge 1$.
\begin{example}[$\psi_{\theta}$-norm of bounded r.v.]\label{eg:boundedr.v.}
    For a r.v. $|X| \le M<\infty$, we have
    \begin{equation*}
        \| X \|_{\psi_{\theta}} = \inf \{ t > 0 : \E e^{|X|^{\theta} / t^{\theta}} \leq 2\} \leq \inf \{ t > 0 : \E e^{M^{\theta} / t^{\theta}} \leq 2\} = M (\log 2)^{-1 / \theta}.
    \end{equation*}
\end{example}

In general, we have following corollary to determine $\|X\|_{\psi_{\theta}}$ based on moment generating functions (MGF). It would be useful for doing statistical inference of $\psi_{\theta}$-norm.
\begin{corollary}\label{gote}
If $\|X\|_{\psi_{\theta}}<\infty$, then $\|X\|_{\psi_{\theta}}=\big( m_{|X|^{\theta}}^{-1}(2) \big)^{-1 / \theta}$ for the MGF $\phi_{Z}(t):=\E e^{tZ}$.
\end{corollary}

\begin{remark}\label{remark1}
    If we observe i.i.d. data $\{X_i\}_{i=1}^{n}$ from a sub-Weibull distribution, one can use the empirical moment generating function (EMGF, \cite{Gbur89}) to estimate the sub-Weibull norm of $X$.  Since the EMGF $\hat{m}_{|X|^{\theta}} (t) = \frac{1}{n} \sum_{i = 1}^n \exp \{ t |X_i|^{\theta}\}$
    converge to MGF $m_{|X|^{\theta}}(t)$ in probability for $t$ in a neighbourhood of zero, the value of the inverse function of EMGF at $2$. Then, under some regularity conditions, $\big( \hat{m}_{|X|^{\theta}} \big)^{-1} (2)$, is a consistent estimate for $\|X\|_{\psi_{\theta}}$.
\end{remark}

    In particular, if we take $\theta = 1$, we get the sub-exponential norm of $X$, which is defined as $\|X\|_{\psi_{1}}=\inf \{t>0: \E \exp (|X| / t) \leq 2\}$. For independent r.v.s $\{X_i\}_{i=1}^{n}$, if $\E X_i = 0$ and $\|X_i\|_{\psi_{1}} < \infty$, by Proposition 4.2 in \cite{Zhang20}, we know $\forall~t \geq 0$
    \begin{equation}\label{exp}
        \pr \bigg(\Big|\sum_{i=1}^{n} X_{i}\Big| \geq t\bigg) \leq 2 \exp \left\{-\frac{1}{4}\left(\frac{t^{2}}{\sum_{i=1}^{n} 2\|X_{i}\|_{\psi_{1}}^{2}} \wedge \frac{t}{\displaystyle\max _{1 \leq i \leq n}\|X_{i}\|_{\psi_{1}}}\right)\right\}.
    \end{equation}

\begin{example} An explicitly calculation of the sub-exponential norm is given in \cite{Gotze2019}, they show that Poisson r.v. $X\sim {\rm{Poisson}}(\lambda)$ has sub-exponential norm $\left\|X \right\|_{\psi_{1}}\le[\log (\log (2) \lambda^{-1}+1)]^{-1}$. And Example \ref{eg:boundedr.v.} with triangle inequality implies
   $$\left\|X-\E X \right\|_{\psi_{1}}\le \left\|X\right\|_{\psi_{1}}+\left\|\E X \right\|_{\psi_{1}} =   \left\|X\right\|_{\psi_{1}}+\frac{\lambda}{{ {\log 2} }} \le[\log (\log (2) \lambda^{-1}+1)]^{-1}+\frac{\lambda}{{ {\log 2} }}$$
based on following useful results.
\end{example}
\begin{proposition}[Lemma A.3 in \cite{Gotze2019}]\label{Gotze2019}
For any $\alpha>0$ and any r.v.s $X, Y$ we have $\|X+Y\|_{\psi_{\theta}}  \leq K_{\alpha}\left(\|X\|_{\psi_{\theta}}+\|Y\|_{\psi_{\theta}}\right)$ and
$$
\|\E X\|_{\psi_{\theta}}  \leq \frac{1}{d_{\alpha}(\log 2)^{1 / \alpha}}\|X\|_{\psi_{\theta}},~~\|X-\E X\|_{\psi_{\theta}}  \leq K_{\alpha}\left(1+\left(d_{\alpha} \log 2\right)^{-1 / \alpha}\right)\|X\|_{\psi_{\theta}},
$$
where $d_{\theta}:=({\theta} e)^{1 / {\theta}} / 2$, $K_{\theta}:=2^{1 /{\theta}}$ if ${\theta}\in(0,1)$ and $K_{{\theta}}=1$ if ${\theta} \geq 1$.
\end{proposition}
To extend Poisson variables, one can also consider concentration for sums of independent heterogeneous negative binomial variables $\{Y_i\}_{i=1}^n$ with probability mass functions:
    \begin{equation}\label{eq:NBR}
        \pr(Y_i = y) = \frac{{\varGamma (y + k_i)}}{{\varGamma (k_i)y!}}{(1 - q_i)^{k_i}} {q_i^{y}}\quad \bigl(q_i \in (0,1),y \in \mathbb{N} \bigr),
    \end{equation}
where $\{k_i \}_{i=1}^n\in (0,\infty)$ are variance-dependence parameters. Here, the mean and variance of $\{Y_i\}_{i=1}^n$ are $\mathrm{{E}}{Y_i}=\frac{{{k_i} q_i}}{{1-q_i}},~\operatorname{{Var}} Y_i=\frac{{{k_i} q_i}}{{{(1-q_i)^2}}}$ respectively. The MGF of $\{Y_i\}_{i=1}^n$ are $\mathrm{E}e^{s Y_i} = \left( \frac{{1 - q_i}}{{1 - q_i{e^{s}}}} \right)^{k_i}$ for $i=1,\cdots,n$. Based on \eqref{exp}, we obtain following results.

\begin{corollary}\label{remark2}
    For any independent r.v.s $\{Y_i\}_{i=1}^n$ satisfying $\|Y_i\|_{\psi_1} < \infty$, $t \geq 0$, and non-random weight $\bfw = (w_1, \cdots, w_n)^{\top}$, we have
    \begin{equation*}
        \pr\bigg( |\sum\limits_{i = 1}^n {{w_i}} ({Y_i} - {\rm{E}}{Y_i})| \ge t\bigg) \leq 2 e^{-\frac{1}{4}\left(\frac{t^{2}}{2 \sum_{i=1}^{n} w_i^2( \|Y_{i}\|_{\psi_{1}} + |\E Y_i / \log 2|)^2} \wedge \frac{t}{\max _{1 \leq i \leq n}|w_i| ( \|Y_{i}\|_{\psi_{1}} + |\E Y_i / \log 2| )}\right)}.
    \end{equation*}
    \[
        \pr\bigg(|\sum\limits_{i = 1}^n {{w_i}} ({Y_i} - {\rm{E}}{Y_i})| > 2{\Big(2t\sum\limits_{i = 1}^n {w_i^2\left\| {{Y_i} - {\rm{E}}{Y_i}} \right\|_{{\psi _1}}^2} \Big)^{1/2}} + 2t\mathop {\max }\limits_{1 \le i \le n} ( {|{w_i}|{{\left\| {{Y_i} - {\rm{E}}{Y_i}} \right\|}_{{\psi _1}}}})\bigg)  \le 2{e^{ - t}}.
    \]
    In particular, if $Y_i$ is independently distributed as $\operatorname{NB}(\mu_i, k_i)$, we have
    \begin{equation}\label{eq:NBRcon}
        \pr \bigg(|\sum_{i=1}^{n} w_i (Y_{i} - \E Y_i) | \geq t \bigg) \leq 2 e^{-\frac{1}{4}(\frac{t^{2}}{2 \sum_{i=1}^{n} w_i^2 a^2(\mu_i, k_i)} \wedge \frac{t}{\max _{1 \leq i \leq n}|w_i|a(\mu_i, k_i) })},
    \end{equation}
    where $a (\mu_i, k_i) := \left[ \log \frac{1 - (1 - q_i) / \sqrt[k_i]{2}}{q_i}\right]^{-1} + \frac{\mu_i}{\log 2}$ with $q_{i}:=\frac{\mu_{i}}{k_{i}+\mu_{i}}.$
\end{corollary}

Corollary \ref{remark2} can play an important role in many non-asymptotic analyses of various estimators. For instance, recently \cite{li2022} uses the above inequality as an essential role for deriving the non-asymptotic behavior of the penalty estimator in the counting data model.

Next, we study moment properties for sub-Weibull random variables. Lemma 1.4 in \cite{Rigollet19} showed that if $X \sim \operatorname{subG}(\sigma^{2})$, then we have: (a). the tail satisfies $\pr (|X|>t) \leq 2e^{-{t^{2}}/{2 \sigma^{2}}}$ for any $t>0$; (b). The (a) implies that moments $\mathrm{E}|X|^{k} \leq(2 \sigma^{2})^{k / 2} k \Gamma(\frac{k}{2})$ and $[k^{ - 1 / 2}(\mathrm{E}(|X|^{k}))^{1 / k}]^2 \leq \sigma^2 e^{2 / e},~k \geq 2$. We extend Lemma 1.4 in \cite{Rigollet19} to sub-Weibull r.v. $X$ satisfying following properties.
\begin{corollary}[Moment properties of sub-Weibull norm]\label{prop:Psub-W}
(a). \text{If}~$\| X\|_{\psi_{\theta}} < \infty$, then $\pr \{ |X| > t \} \le 2 e^{-(t / \|X\|_{\psi_{\theta}})^{\theta}}$ for all $t \ge 0$; and then  $\E |X|^k \le 2 \|X\|_{\psi_{\theta}}^{k} \Gamma(\frac{k}{\theta}+1)$ for all $k\ge 1$. (2). Let ${C_\theta } := \mathop {\max }\limits_{k \ge 1} {\left( {\frac{{2\sqrt {2\pi } }}{\theta }} \right)^{1/k}}{\left( {\frac{k}{\theta }} \right)^{1/(2k)}}$, for all $k \geq 1$ we have $(\E |X|^k)^{1/k} \le  {C_\theta }{( {\theta {e^{11/12}}} )^{ - 1/\theta }} \| X\|_{\psi_{\theta}}{k^{1/\theta }}.$
\end{corollary}

Particularly, sub-Weibull r.v.s reduce to sub-exponential or sub-Gaussian r.v.s when $\theta$ = 1 or 2. It is obvious that the smaller $\theta$ is, the heavier tail the r.v. has. A r.v. is called heavy-tailed if its distribution function fails to be bounded by a decreasing exponential function, i.e. 
\begin{center}
$\int e^{\lambda x} d F(x)=\infty, \forall \lambda>0$ (the tail decays slower than some exponential r.v.s);
\end{center}
see \cite{Foss2011}. Hence for sub-Weibull r.v.s, we usually focus on the the sub-Weibull index $\theta \in (0, 1)$. A simple example that the heavy-tailed distributions arises when we work more production on sub-Gaussian r.v.s. Via a power transform of $|X|$, the next corollary explains the relation of sub-Weibull norm with parameter $\theta$ and $r\theta$, which is similar to Lemmas 2.7.6 of \cite{Vershynin18} for sub-exponential norm.

\begin{corollary}\label{lem:2theta}
    For any $\theta,r \in(0, \infty),$ if  $X \sim \operatorname{subW}(\theta)$, then $|X|^r \sim \operatorname{subW}(\theta/r)$. Moreover,
    \begin{equation}\label{eq:2theta}
         \left\||X|^{r}\right\|_{\psi_{\theta/r}}=\|X\|^r_{\psi_{\theta}}.
    \end{equation}
    Conversely, if  $X \sim \operatorname{subW}(r\theta)$, then $X^r \sim \operatorname{subW}(\theta)$ with $\left\|X^{r}\right\|_{\psi_{\theta}} =\|X\|_{\psi_{r\theta}}^{r}$.
\end{corollary}

 By Corollary \ref{lem:2theta}, we obtain that $d$-th root of the absolute value of sub-Gaussian is $\operatorname{subW}(2d)$ by letting $r=1/d$. Corollary \ref{lem:2theta} can be extended to product of r.v.s, from Proposition D.2 in \cite{Kuchibhotla18} with the equality replacing by inequality, we state it as the following proposition.

\begin{proposition}\label{co: dproduct subw}
 If $\{W_{i}\}_{i=1}^d$ are (possibly dependent) r.vs satisfying $\left\|W_{i}\right\|_{\psi_{\alpha_{i}}}<$ $\infty$ for some $\alpha_{i}>0,$ then
$$
\|\prod_{i=1}^{d} W_{i}\|_{\psi_{\beta}} \leq \prod_{i=1}^{d}\|W_{i}\|_{\psi_{\alpha_{i}}} \text { where } \frac{1}{\beta}:=\sum_{i=1}^{d} \frac{1}{\alpha_{i}}.
$$
\end{proposition}

For multi-armed bandit problems in reinforcement learning, \cite{Hao19} move beyond sub-Gaussianity and consider the reward under sub-Weibull distribution which has a much weaker tail. The corresponding concentration inequality (Theorem 3.1 in \cite{Hao19}) for the sum of independent sub-Weibull r.v.s is illustrated as follows.
\begin{proposition}[Concentration inequality for sub-Weibull distribution]\label{thm:orlicz_concentra}
	Suppose  $\{X_i\}_{i=1}^n$ are independent sub-Weibull random variables with $\|X_i- \E X_i\|_{\psi_{\theta}}\leq v$. Then there exists absolute constants $C_{1\theta}$ and $C_{2\theta}$ only depending on ${\theta}$ such that with probability at least $1-{e^{ - t}}$:
\[\left|\frac{1}{n}\sum\limits_{i = 1}^n \frac{{X_i} - \E{X_i}}{v} \right| \le {C_{1\theta }}{\left( {\frac{t}{n}} \right)^{1/2}} + {C_{2\theta }}{\left( {\frac{t}{n}} \right)^{1/\theta }} = \left\{ \begin{array}{l}
O({n^{ - 1/\theta }}),\theta > 2\\
O({n^{ - 1/2}} ),0 < \theta  \le 2
\end{array} \right..\]
\end{proposition}
The weakness in the Proposition \ref{thm:orlicz_concentra} is that the upper bound of $S_n^{\bm{a}}:=\sum_{i=1}^na_i Y_i-\mathrm E(\sum_{i=1}^na_i Y_i)$ is up to a unknown constants $C_{1\theta}, C_{2\theta}$. In the next section, we will give a constants-specified and high probability upper bound for $|S_n^{\bm{a}}|$, which improve Proposition \ref{thm:orlicz_concentra} and  is sharper than Theorem 3.1 in \cite{Kuchibhotla18}.

\subsection{Main results: concentrations for sub-Weibull summation}

Based on the exponential moment condition, the Chernoff's tricks implies the following sub-exponential concentrations from Proposition 4.2 in \cite{Zhang20}.
\begin{proposition}\label{remark3}
    For any independent r.v.s $\{Y_i\}_{i=1}^n$ satisfying $\|Y_i\|_{\psi_1} < \infty$, $t \geq 0$, and non-random weight $\bfw = (w_1, \cdots, w_n)^{\top}$, we have
\begin{center}
$\pr(|\sum\limits_{i = 1}^n {{w_i}} ({Y_i} - {\rm{E}}{Y_i})| > 2{(2t\sum\limits_{i = 1}^n {w_i^2\left\| {{Y_i} - {\rm{E}}{Y_i}} \right\|_{{\psi _1}}^2} )^{1/2}} + 2t\mathop {\max }\limits_{1 \le i \le n} ( {|{w_i}|{{\left\| {{Y_i} - {\rm{E}}{Y_i}} \right\|}_{{\psi _1}}}}))  \le 2{e^{ - t}}.$
\end{center}
\end{proposition}
But it is not easy to extend to sub-Weibull distributions. From Corollary \ref{lem:2theta}, ${Y_i} \sim \operatorname{subW}(\theta) \Rightarrow  |{Y_i}|^{1/\theta} \sim \operatorname{subW}(1)$. The MGF of $|{Y_i}|^{1/\theta}$ satisfies $\mathrm{E}e^{ \lambda^{1/\theta} |{Y_i}|^{1/\theta} }\le e^{\lambda^{1/\theta} K^{1/\theta} },~|\lambda| \le \frac1{K}$ for some constant ${K}>0$. The bound of $\mathrm{E}e^{ \lambda^{1/\theta} |{Y_i}|^{1/\theta} }$ with $\theta \neq 1$ or 2 is not directly applicable for deriving the concentration of $\sum_{i = 1}^n {{w_i}} ({Y_i} - {\rm{E}}{Y_i})$ by using the independence and Chernoff's tricks, since the MGF of Weibull r.v. do not has closed form as exponential function. Thanks to the tail probability derived by Orlicz-type norms, instead of using the upper bound for MGF, an alternative method is given by \cite{Kuchibhotla18} who defines the so-called Generalized Bernstein-Orlicz (GBO) norm.  And the GBO norm can help us to derive tail behaviours for sub-Weibull r.v.s.

\begin{definition}[GBO norm]\label{def:GBOnorm}
Fix $\alpha > 0$ and $L \ge 0$. Define the function $\Psi_{\theta, L}(\cdot)$ as the inverse function $\Psi_{\theta, L}^{-1}(t) := \sqrt{\log (t+1)} + L\left(\log (t+1)\right)^{1/\theta}~\mbox{for all}~t\ge 0.$ The GBO norm of a r.v. $X$ is then given by
$\|X\|_{\Psi_{\theta, L}}:=\inf \{\eta>0: \mathrm{E}[{\Psi_{\theta, L}}(|X| / \eta)] \leq 1\}.$
\end{definition}
The monotone function $\Psi_{\theta, L}(\cdot)$ is motivated by the classical Bernstein's inequality for sub-exponential r.v.s. Like the sub-Weibull norm properties Corollary \ref{prop:Psub-W}, the following proposition in \cite{Kuchibhotla18} allows us to get the concentration inequality for r.v. with finite GBO norm.
\begin{proposition}\label{prop:GBO norm}
If $\|X\|_{\Psi_{\theta, L}}< \infty$, then ${\pr}(|X| \ge \|X\|_{\Psi_{\theta, L}}\{\sqrt{t} + Lt^{1/\theta}\}) \le 2e^{-t}~\forall~t\ge 0.$
\end{proposition}

With an upper bound of GBO norm, we could easily derive the concentration inequality for a single sub-Weibull r.v. or even the sum of independent sub-Weibull r.v.s. The sharper upper bounds for the GBO norm is obtained for the sub-Weibull summation, which refines the constant in the sub-Weibull concentration inequality. Let $||X||_p:=(\E |X|^p)^{1/p}$ for all integer $p \ge 1$. First, by truncating more precisely, we obtain a sharper upper bound for $||X||_p$, comparing to Proposition C.1 in \cite{Kuchibhotla18}.
\begin{corollary}\label{lem2}
    If $\|X\|_p \leq C_1 \sqrt{p} + C_2 p^{1/\theta}$ for $p \geq 2$ and constants $C_1$, $C_2$, then
    \begin{equation*}
        \|X\|_{\Psi_{\theta, K}} \leq \gamma e C_1
    \end{equation*}
    where $K = \gamma^{2 / \theta} C_2  /  (\gamma C_1)$ and $\gamma \approx 1.78$ is the minimal solution of 
    \begin{equation*}
        \left\{ k > 1:{e^{{2k^{ - 2}}}} - 1 + \frac{{{e^{2(1 - {k^2})/{k^2}}}}}{{{k^2} - 1}} \le 1\right\}.
    \end{equation*}
\end{corollary}

The proof can be seen in the Appendix. In below, we need the moment estimation for sums of independent symmetric r.v.s.

\begin{lemma}[Khinchin-Kahane Inequality, Theorem 1.3.1 of \cite{Pena12}]\label{lemKK}
   Let $\left\{a_{i}\right\}_{i=1}^{n}$ be a finite non-random sequence, $\left\{\varepsilon_{i}\right\}_{i=1}^{n}$ be a sequence of independent Rademacher variables and $1<p<q<\infty .$ Then
$
\left\|\sum_{i=1}^{n} \varepsilon_{i} a_{i}\right\|_{q} \leq\left(\frac{q-1}{p-1}\right)^{1 / 2}\left\|\sum_{i=1}^{n} \varepsilon_{i} a_{i}\right\|_{p}.
$
\end{lemma}

\begin{lemma}[Theorem 2 of \cite{Latala97}]\label{lem3}
    Let $\left\{X_{i}\right\}_{i=1}^{n}$ be a sequence of independent \emph{symmetric} r.v.s, and $p \geq 2$. Then, $\frac{e-1}{2 e^{2}}\left\|\left(X_{i}\right)\right\|_{p} \leq\left\|X_{1}+\cdots+X_{n}\right\|_{p} \leq e\left\|\left(X_{i}\right)\right\|_{p},$ where $ \left\|\left(X_{i}\right)\right\|_{p} :=\inf \{t>0: \sum_{i = 1}^n \log \phi_{p}\left({X_{i}}/{t}\right) \leq p\}$ with $\phi_p(X) := \mathrm{E} |1 + X|^p.$
\end{lemma}

\begin{lemma}[Example 3.2 and 3.3 of \cite{Latala97}]\label{lem4}
    Assume $X$ be a symmetric r.v. satisfying $\pr \left( |X| \geq t\right) =e^{ - N(t)}$. For any $t \geq 0$, we have
    \begin{itemize}
        \item[\rm{(a)}] If $N(t)$ is concave, then $\log \phi_{p}(e^{-2} t X) \leq p M_{p, X}(t) := (t^{p}\|X\|_{p}^{p}) \vee (p t^{2}\|X\|_{2}^{2})$.

        \item[\rm{(b)}] For convex $N(t)$, denote the convex conjugate function $N^{*}(t) := \sup_{ s>0} \{t s-N(s)\}$ and $M_{p, X}(t)=\left\{\begin{array}{ll}p^{-1} N^{*}(p|t|), & \text { if  } p|t| \geq 2 \\ p t^{2}, & \text { if  } p|t|<2.\end{array}\right.$ Then  $\log \phi_{p}(t X / 4) \leq p M_{p, X}(t)$.
    \end{itemize}
\end{lemma}
With the help of three lemmas above, we can obtain the main results concerning the shaper and constant-specified concentration inequality for the sum of independent sub-Weibull r.v.s.
\begin{theorem}[Concentration for sub-Weibull summation]\label{thm:SumNewOrliczVex}
Let $\gamma$ be given in Corollary \ref{lem2}. If $\left\{X_{i}\right\}_{i=1}^{n}$ are independent centralized r.v.s 
such that $\|X_i\|_{\psi_{\theta}} < \infty$ for all $1\le i\le n$ and some $\theta > 0$, then for any weight vector $\bm w= (w_1, \ldots, w_n)\in\mathbb{R}^n$, the following bounds holds true:
    \begin{itemize}
        \item[\rm{(a)}] The estimate for GBO norm of the summation:
    \begin{center}
    $\left\|\sum_{i=1}^n w_iX_i\right\|_{\Psi_{\theta, L_n(\theta,\bm{b}_X)}} \le \gamma eC(\theta)\| \bm{b}_X\|_2$,
    \end{center}
    where ${\bm b}_X= (w_1\|X_1\|_{\psi_{\theta}}, \ldots, w_n\|X_n\|_{\psi_{\theta}})^{\top}\in\mathbb{R}^n$, with
\begin{equation*}
    C(\theta) :=  \begin{cases}
        2 \left[ \log^{1 / \theta} 2 + e^3 \left( \Gamma^{1 / 2} \left( \frac{2}{\theta} + 1\right) + 3^{\frac{2 - \theta}{3 \theta}} \sup_{p \geq 2} p^{ - \frac{1}{\theta}} \Gamma^{1 / p} \left( \frac{p}{\theta} + 1\right) \right) \right], &\mbox{ if } \theta \leq 1,\\
         2[4e + (\log 2)^{1/\theta}], &\mbox{ if } \theta > 1;
    \end{cases}
\end{equation*}
and $L_n (\theta, \bm b) = \gamma^{2 / \theta} A(\theta) \frac{\| \bm b\|_{\infty}}{\| \bm b\|_2} 1 \{ 0 < \theta \leq 1\} + \gamma^{2 / \theta} B(\theta) \frac{\| \bm b\|_{\beta}}{\| \bm b\|_2} 1 \{ \theta > 1\}$ where $B(\theta)=:\frac{2 e \theta^{- 1 / \theta} \left( 1 - \theta^{-1}\right)^{1 / \beta} }{4 e + (\log 2)^{1 / \theta}} $ and $A(\theta)=: \inf\limits_{p \geq 2} \frac{e^{3} 3^{\frac{2 - \theta}{3 \theta}}p^{- 1 / \theta} \Gamma^{1 / p} \left( \frac{p}{\theta} + 1\right)  }{2 [ \log^{1 / \theta} 2 + e^3( \Gamma^{1 / 2}( \frac{2}{\theta} + 1) + 3^{\frac{2 - \theta}{3 \theta}} \sup_{p \geq 2} p^{- {1}/{\theta}} \Gamma^{1 / p} ( \frac{p}{\theta} + 1) ) ]} $. For the case $\theta > 1$, $\beta$ is the H\"{o}lder conjugate satisfying $1/\theta + 1/\beta = 1$.

\item[\rm{(b)}] Concentration for sub-Weibull summation:
\begin{equation}\label{eq:sWc}
    \pr \bigg(|\sum\limits_{i = 1}^n {{w_i}{X_i}} | \ge 2eC(\theta ){\| \bm{b}_X \|_2}\{ \sqrt t  + {L_n}(\theta, \bm{b}_X){t^{1/\theta }}\} \bigg) \le 2{e^{ - t}}.
\end{equation}

\item[\rm{(c)}] Another form of for $\theta \ne  2$:
\begin{align*}
    \pr \bigg(  {| {\sum\limits_{i = 1}^n {{w_i}} {X_i}} | \ge s} \bigg) &\le 2\exp \left\{ { - \bigg( {\frac{{{s^\theta }}}{{{{\big[4eC(\theta ){{\left\| \bm{b}_X \right\|}_2}{L_n}(\theta, \bm{b}_X)\big]}^\theta }}} \wedge \frac{{{s^2}}}{{16{e^2}{C^2}(\theta )\left\| \bm{b}_X \right\|_2^2}}} \bigg)} \right\}\\
    (\theta  < 2) & = \left\{ {\begin{array}{*{20}{l}}
    {2{e ^{ - {s^2}/16{e^2}{C^2}(\theta )\left\| b \right\|_2^2}}, \text{ if } s \le 4eC(\theta ){{\left\| \bm{b}_X \right\|}_2}L_n^{\theta /(\theta -  2)}(\theta, \bm{b}_X)}\\
    {2{e ^{ - {s^\theta }/{{[4eC(\theta ){{\left\| \bm{b}_X \right\|}_2}{L_n}(\theta, \bm{b}_X)]}^\theta }}}, \text{ if } s > 4eC(\theta ){{\left\| \bm{b}_X \right\|}_2}L_n^{\theta /(\theta -  2)}(\theta, \bm{b}_X)};
    \end{array}} \right.\\
    (\theta  > 2) &= \left\{ {\begin{array}{*{20}{l}}
    {2{e ^{ - {s^\theta }/{{[4eC(\theta ){{\left\| \bm{b}_X \right\|}_2}{L_n}(\theta, \bm{b}_X)]}^\theta }}}, \text{ if } s < 4eC(\theta ){{\left\| \bm{b}_X \right\|}_2}L_n^{\theta /(2-\theta )}(\theta, \bm{b}_X)}\\
    {2{e ^{ - {s^2}/16{e^2}{C^2}(\theta )\left\| \bm{b}_X \right\|_2^2}}, \text{ if } s \ge 4eC(\theta ){{\left\| \bm{b}_X \right\|}_2}L_n^{\theta /(2-\theta )}(\theta, \bm{b}_X)}.
    \end{array}} \right.
\end{align*}
\end{itemize}
\end{theorem}
\begin{remark}
 The constant $C(\theta)$ in Theorem \ref{thm:SumNewOrliczVex} can be improved as $C(\theta)/2$ under symmetric assumption of sub-Weibull r.v.s  $\{X_i\}_{i=1}^n$. Moreover, by the \textit{improved symmetrization theorem} (Theorem 3.4 in \cite{Kashlak2018}), one can replace the constant $C(\theta)$ in Theorem \ref{thm:SumNewOrliczVex} by a sharper constant $(1+o(1))C(\theta)/2$.  Theorem \ref{thm:SumNewOrliczVex} (b) also implies a potential empirical upper bound for $\sum_{i = 1}^n w_i X_i$ for independent sub-Weibull r.v.s $\{X_i\}_{i=1}^n$, because the only unknown variable in $2 e C(\theta) \| \bm{b}_X\|_2 \{ \sqrt{t} + L_n(\theta) t^{1 / \theta}\}$ is $\bm{b}_X$. From Remark \ref{remark1}, estimating ${\bm b}_X$ is possible for i.i.d. observation $\{X_i\}_{i=1}^n$.
\end{remark}
 \begin{remark}
 Compared with the newest result in \cite{Kuchibhotla18}, our method do not use the crude String's approximation will give sharper concentration. For example, suppose $X_1, \ldots, X_{10}$ are i.i.d. r.v.s with mean $\mu$ and $\| X_1 - \mu \|_{\psi_{\theta}} = 1$. Here we set $\theta = 0.5$, $X$ is heavy-tailed (for example set the density of $X$ as $f(x) = \frac{1}{2 \sqrt{x}} e^{- \sqrt{x}}\cdot 1 (x \geq 0)$). We find that $C(\theta) \approx 2825.89$, $A(\theta) \approx 0.07$, and $L_{10}(\theta, \mathbf{1}_{10}^{\top}) = 0.23$. Hence, $95\%$ confidence interval in our method will be
\[
    \mu \in \overline{X} \pm 2 e \times 2118.80,
\]
while the 95\% confidence interval in Theorem 3.1 of \cite{Kuchibhotla18} is evaluated as
\[
    \mu \in \overline{X} \pm 2e \times 3969.94.
\]
In this example, it can be seen that our method does give a much better (tighter) confidence interval.
\end{remark}

\begin{remark}
    Theorem \ref{thm:SumNewOrliczVex} (b) generalizes the sub-Gaussian concentration inequalities, sub-exponential concentration inequalities, and Bernstein's concentration inequalities with Bernstein's moment condition. For $\theta < 2$ in Theorem \ref{thm:SumNewOrliczVex} (c), the tail behaviour of the sum is akin to a sub-Gaussian tail for small $t$, and the tail resembles the exponential tail for large $t$; {{For $\theta  > 2$, the tail behaves like a Weibull r.v. with tail parameter $\theta$ and the tail of sums match that of the sub-Gaussian tail for large $t$}}. The intuition is that the sum will concentrate around zero by the Law of Large Number. Theorem \ref{thm:SumNewOrliczVex} shows that the convergence rate will be faster for small deviations from the mean and will be slower for large deviations from the mean.
\end{remark}

\begin{remark}
    Recently, similar result presented in \cite{vladimirova20} is that
    \begin{equation*}
        \pr \bigg( \Big| \sum_{i = 1}^n X_i \Big| > x \bigg) \leq \exp \bigg\{ - \Big( \frac{x}{n K_{\theta}}\Big)^{1 / \theta} \bigg\},~\text{for}~x\ge n K_{\theta}
    \end{equation*}
    where $K_{\theta}$ is some constants only depends on $X$ and $\theta$ ($K_{\theta}$ can be obtained by Proposition \ref{prop:Psub-W}). But it is obvious to see this large derivation result cannot guarantee a $\sqrt{n}$-convergence rate (as presented in Proposition \ref{thm:orlicz_concentra}) whereas our result always give a $\sqrt{n}$-convergence rate, as presented in Theorem \ref{thm:SumNewOrliczVex} (c) and Proposition \ref{thm:orlicz_concentra}.
\end{remark}

\subsection{Sub-Weibull parameter}
In this part, a new \emph{sub-Weibull parameters} is proposed, which is enable of recovering the tight concentration inequality for single non-zero mean random vector. Similar to characterizations of sub-Gaussian r.vs. in Proposition 2.5.2 of \cite{Vershynin18}, sub-Weibull r.vs. has the equivalent definitions.
\begin{proposition}[Characterizations of sub-Weibull r.v., \cite{Wong17}]\label{th:subWeibull}
  Let $X$ be a r.v., then the following properties are equivalent. (1). The tails of $X$ satisfy $\pr(|X| \ge x) \le e^{ - (x/ K_1)^{\theta}},~\text{for all } x \ge 0$; (2).  The moments of $X$ satisfy $\|X\|_k:=(\mathrm{E}|X|^{k})^{1 / k} \leq K_{2} k^{1 /\theta}~\text{for all } k \ge 1 \wedge \theta$; (3).  The MGF of $|X|^{1/\theta}$ satisfies $\mathrm{E}e^{ \lambda^{1/\theta} |X|^{1/\theta} }\le e^{\lambda^{1/\theta} K_3^{1/\theta} }$ for $|\lambda| \le \frac1{K_3}$; (4). ${\rm{E}}e^ {|X/K_4|^{1/\theta}}  \le 2.$
\end{proposition}

From the upper bound of $(\E |X|^k)^{1/k}$ in Proposition \ref{th:subWeibull}(2), an alternative definition of the sub-Weibull norm
$\|X\|_{\psi_{\theta}}:=\sup _{k \geq 1}k^{-1 /\theta}(\mathrm{E}|X|^{k})^{1 / k} $
  is given by \cite{Wong17}. Let ${\theta}=1$. An alternative definition of the sub-exponential norm is
$\|X\|_{\psi_1} := \sup_{k \ge 1} k^{-1} ({\rm{E}} |X|^k)^{1/k}$
see Proposition 2.7.1 of \cite{Vershynin18}. The sub-exponential r.v. $X$ satisfies equivalent properties in Proposition~\ref{th:subWeibull} (Characterizations of sub-exponential with ${\theta}=1$). However, these definition is not enough to obtain the sharp parameter as presented in the sub-Gaussian case. Here, we redefine the sub-Weibull parameter by our Corollary \ref{prop:Psub-W}(a).

\begin{definition}[Sub-Weibull r.v.,$X  \sim  \operatorname{subW}(\theta,v)$]\label{def.SW}
Define the sub-Weibull norm
\begin{center}
${\left\| X \right\|_{{\varphi _\theta }}} = {\sup }_{k \ge 1} {\left( {{\rm{E}}|X|^{\theta k}}/k! \right)^{1/(\theta k)}}.$
\end{center}
\end{definition}
We denote the sub-Weibull r.v. as $X  \sim  \operatorname{subW}(\theta,v)$ if $v={\left\| X \right\|_{{\varphi _\theta }}}<\infty$ for a given $\theta>0$. For $\theta \ge 1$, the ${\left\| \cdot \right\|_{{\varphi _\theta }}}$ is a norm which satisfies triangle inequality by Minkowski's inequality: $\mathrm{E}(|X+Y|^{r})^{1 / r} \le [\mathrm{E}(|X|^{r})]^{1 / r}+[\mathrm{E}(|Y|^{r})]^{1 / r}$,$(r \ge 1)$ comparing to Proposition \ref{Gotze2019}. Definition \ref{def.SW} is free of bounding MGF, and it avoids Stirling's approximation in the proof of the tail inequality. We obtain following main results for this moment-based norm.
\begin{corollary}\label{cor:type1SW}
If ${\left\| X \right\|_{{\varphi _\theta }}}< \infty$, then
$\pr \{ |X| > t \} \le 2\exp \{  - \frac{{{t^\theta }}}{{2\left\| X \right\|_{{\varphi _\theta }}^\theta }}\} ~\text{for all } t \ge 0.$
\end{corollary}

\begin{theorem}[sub-Weibull concentration]\label{thm:sub-con}
Suppose that there are $n$ independent sub-Weibull r.v.s $X_{i} \sim  \operatorname{subW}(\theta,v_i)$ for $i=1,2,\cdots,n$. We have
\[\pr \left( \left| \sum\limits_{i = 1}^n {{X_i}}  \right| \ge t \right) \le \exp \left\{  - \frac{{\theta {e^{11/12}}{t^\theta }}}{{2{[e(\sum\limits_{i = 1}^n {{v_i}} ){C_\theta }]^\theta }}}\right\},~\text{ for } t \ge e(\sum\limits_{i = 1}^n {{v_i}} ){C_\theta }{{({2^{ - 1}}\theta {e^{11/12}})}^{ - 1/\theta }},\]
and $\pr \left( {\left| {\frac{1}{n}\sum\limits_{i = 1}^n {{X_i}} } \right| \le e\bar v{2^{1/\theta }}{C_\theta }{{\left( {\frac{{\log ({\alpha ^{ - 1}})}}{{\theta {e^{11/12}}}}} \right)}^{1/\theta }}} \right) \ge 1 - \alpha  \in (1 - {e^{ - 1}},1]$. Moreover, we have
\[\pr \left(| \sum\limits_{i = 1}^n {{X_i}} | \ge e{( {\sum\limits_{i = 1}^n ({\rm{E}}|{X_i}|)^t })^{1/t}} + e(\sum\limits_{i = 1}^n {{v_i}} ){2^{1/\theta }}{C_\theta }{( {\frac{t}{{\theta {e^{11/12}}}}} )^{1/\theta }}\right)\le e^{-t},~\forall~t \ge 0.\]
\end{theorem}

The proof of Theorem \ref{thm:sub-con} can be seen in section \ref{thm2.proof}. The concentration in this Theorem \ref{thm:sub-con} will serve a critical role in many statistical and machine learning literature. For instance, the sub-Weibull concentrations in \cite{Hao19} contain unknown parameters, which makes the algorithm for general sub-Weibull random rewards is infeasible. However, when using our results, it will become feasible as we give explicit constants in these concentrations.

Importantly, the sub-exponential parameter is a special case of sub-Weibull norm by letting $\theta=1$. Denote the \textbf{sub-exponential parameter} for r.v $X$ as
\begin{center}
${\left\| X \right\|_{{\varphi _1}}} :=\sup\limits_{k \ge 1} {\left( {\frac{{{\rm{E}}|X{|^k}}}{{k!}}} \right)^{1/k}}$.
\end{center}

We denote $X  \sim  \operatorname{sE}_{\varphi_1}(v)$ if $v=\|X\|_{\varphi_2}$. For exponential r.v. $X \sim {\rm{Exp}}(\mu)$, the moment is ${{\rm{E}}{X^k}} = {{k!}}{{{\lambda ^k}}}$ and ${\left\| X \right\|_{{\varphi _1}}}=\lambda$. Another case of sub-Weibull norm is $\theta=2$, which defines \textbf{sub-Gaussian parameter}:
\begin{center}
${\left\| X\right\|_{{\varphi _2}}} :=\sup\limits_{k \ge 1} {\left( {\frac{{{\rm{E}}|X|^{2k}}}{{k!}}} \right)^{1/{2k}}}\ge (\operatorname{Var}X)^{1/2}$.
\end{center}

Like the generalized method of moments, we can give the higher-moment estimation procedure for the norm $\| X \|_{\varphi_2}$. Unfortunately, the method in Remark \ref{remark1} for estimating MGF is not stable in the simulation since the exponential function has a massive variance in some cases.
\begin{itemize}
\item \textbf{Estimation procedure for $\| X \|_{\varphi_2}$ and $\| X \|_{\varphi_1}$}. \emph{Consider
\begin{equation}\label{est3}
    \widehat{\| X \|}_{\varphi_2}= \sup_{k \ge 1}  \Big( \frac{1}{n \times k!} \sum_{i = 1}^n |X_i|^{2k}\Big)^{1 / (2k)}, \widehat{\|X\|}_{\varphi_{1}}=\sup _{k \geq 2}\left(\frac{1}{k !} \cdot \frac{1}{n} \sum_{i=1}^{n}\left|X_{i}\right|^{k}\right)^{1 / k}
\end{equation}
as a discrete optimization problem. We can take $k_{\max}$ big enough to minimize
\begin{center}
$\Big( \frac{1}{n \times k!} \sum_{i = 1}^n |X_i|^{2k}\Big)^{1 / (2k)},~\left(\frac{1}{k !} \cdot \frac{1}{n} \sum_{i=1}^{n}\left|X_{i}\right|^{k}\right)^{1 / k}$ on $k \in \{1, \ldots, p_{\max} \}$.
\end{center}}
\end{itemize}
At the first glimpse, the bigger $p$ is, the larger $n$ is required in this method. Nonetheless, often, most of  common distributions only require a median-size of $p$ to give a relatively good result, then only the median-size of $n$ in turn is required.
For standard Gaussian random, centralized Bernoulli (successful probability $\mu = 0.3$), and uniform distributed (on $[-1, 1]$) variable $X$,
\begin{equation*}
\| X \|_{\varphi_2}= \sqrt{2} \left[ \frac{\Gamma \big( (1 + p) / 2\big) }{\Gamma (1/2) \Gamma (1 + p / 2)}\right]^{1 / p}, \quad \left[ \frac{\mu(1 - \mu)^p + (1 - \mu) \mu^p}{\Gamma (p / 2 + 1)} \right]^{1 / p}, \quad \frac{\Gamma^{-1 / p}({p} / {2} + 1)}{ (p + 1)^{1 / p}}.
\end{equation*}

\begin{figure}[!htb]
    \minipage{0.333\textwidth}
        \includegraphics[width=\linewidth]{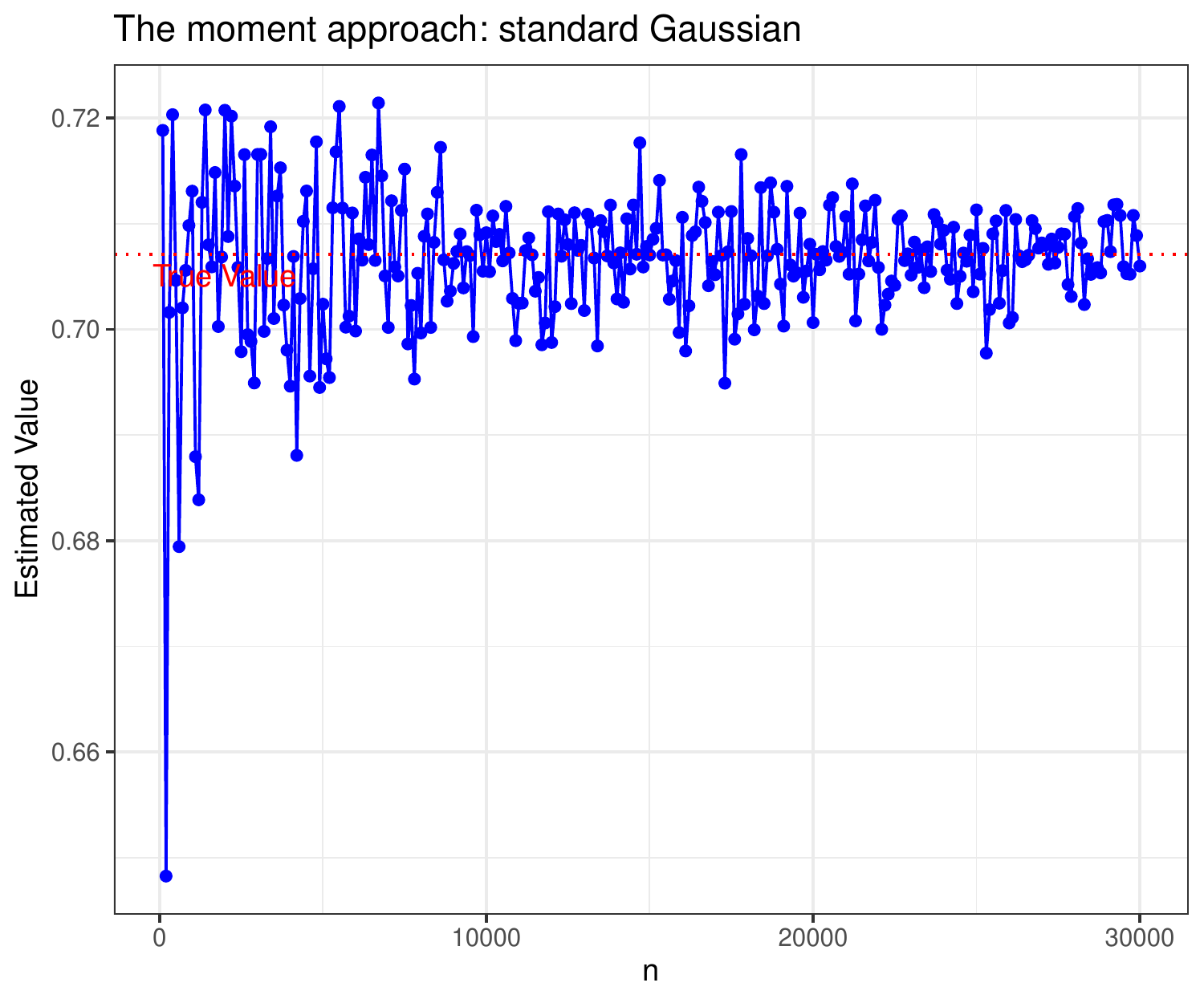}
    \caption{standard Gaussian}
    \label{figure1}
    \endminipage\hfill
    \minipage{0.333\textwidth}
      \includegraphics[width=\linewidth]{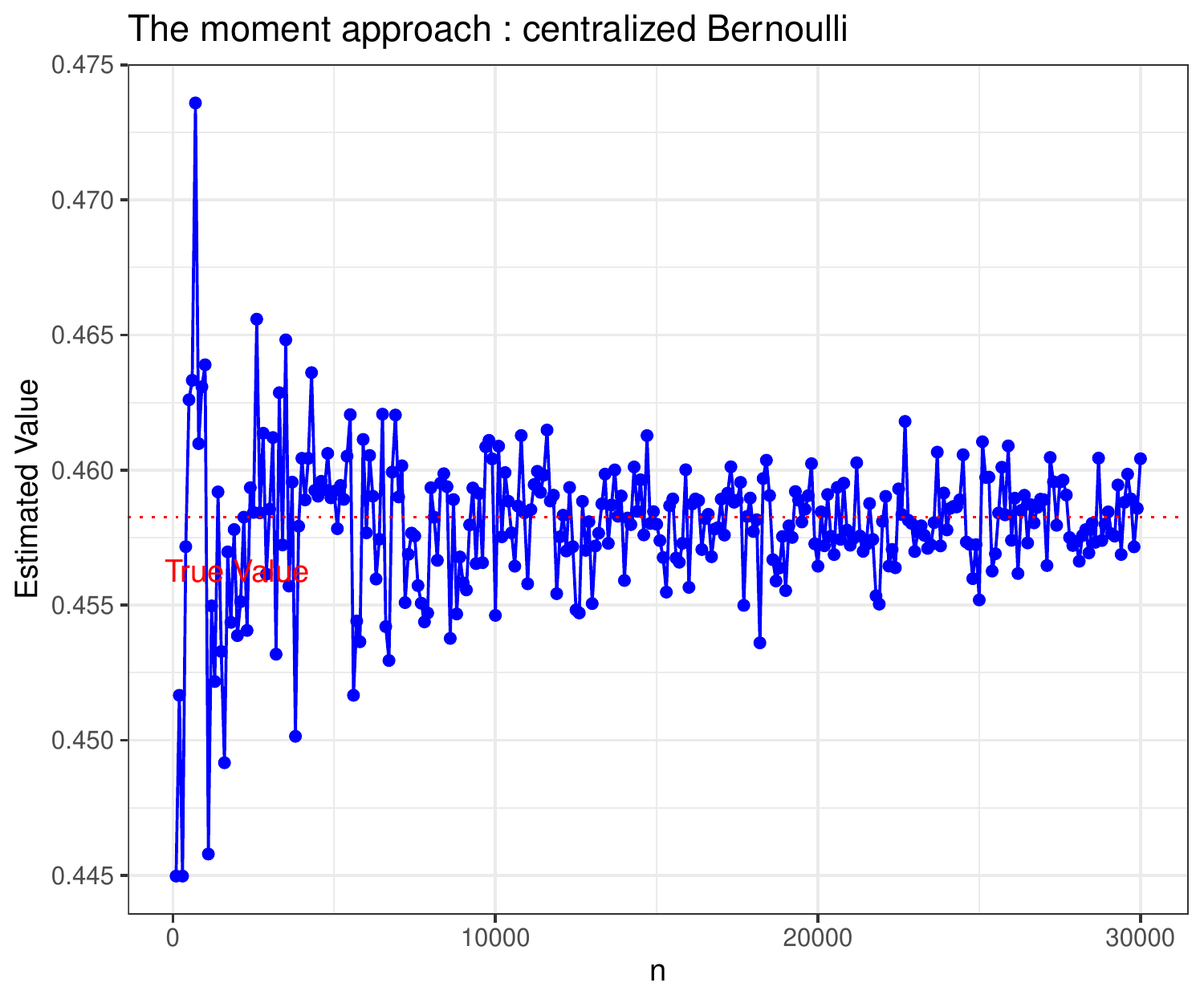}
     \caption{centralized Bernoulli}
     \label{figure2}
    \endminipage\hfill
    \minipage{0.333\textwidth}
      \includegraphics[width=\linewidth]{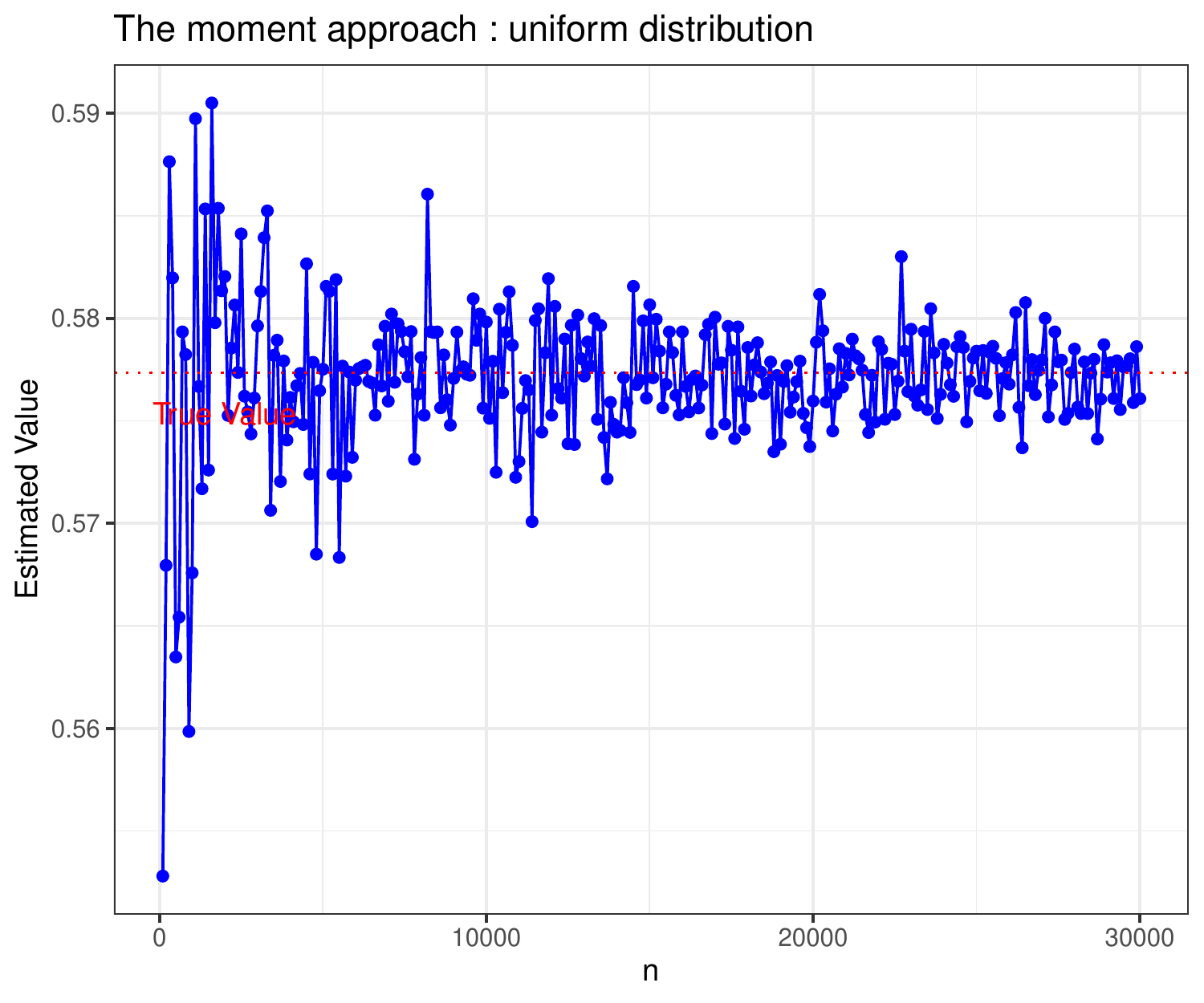}
      \caption{Uniform on $[-1, 1]$}
      \label{figure3}
    \endminipage
\end{figure}

   It can be shown that $\| X\|_{\varphi_2}  \approx 1, 0.4582576, 0.5773503.$ The Figure \ref{figure1}, Figure \ref{figure2}, and Figure \ref{figure3}  show the estimated value from different $n$ under estimate method (\ref{est3}) for the three distributions mentioned above. The estimate method (\ref{est3}) is a correct estimated method for sub-Gaussian parameter to our best knowledge.

\begin{figure}[!htb]
    \minipage{0.5\textwidth}
    \begin{center}
        \includegraphics[width=.7\linewidth]{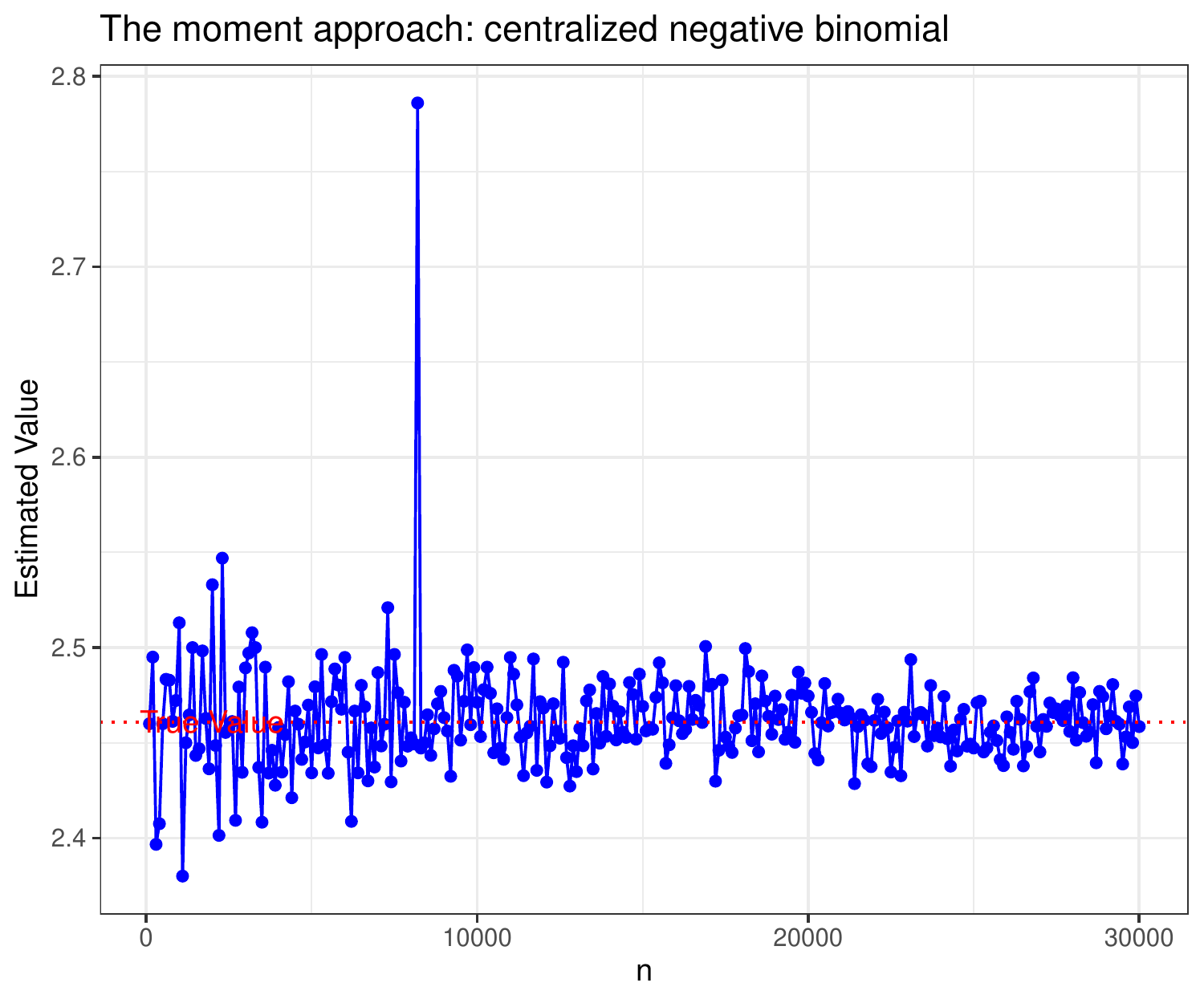}
        \caption{centralized negative binomial}
        \label{figure4}
    \end{center}
    \endminipage
    \minipage{0.5\textwidth}
    \begin{center}
        \includegraphics[width=.7\linewidth]{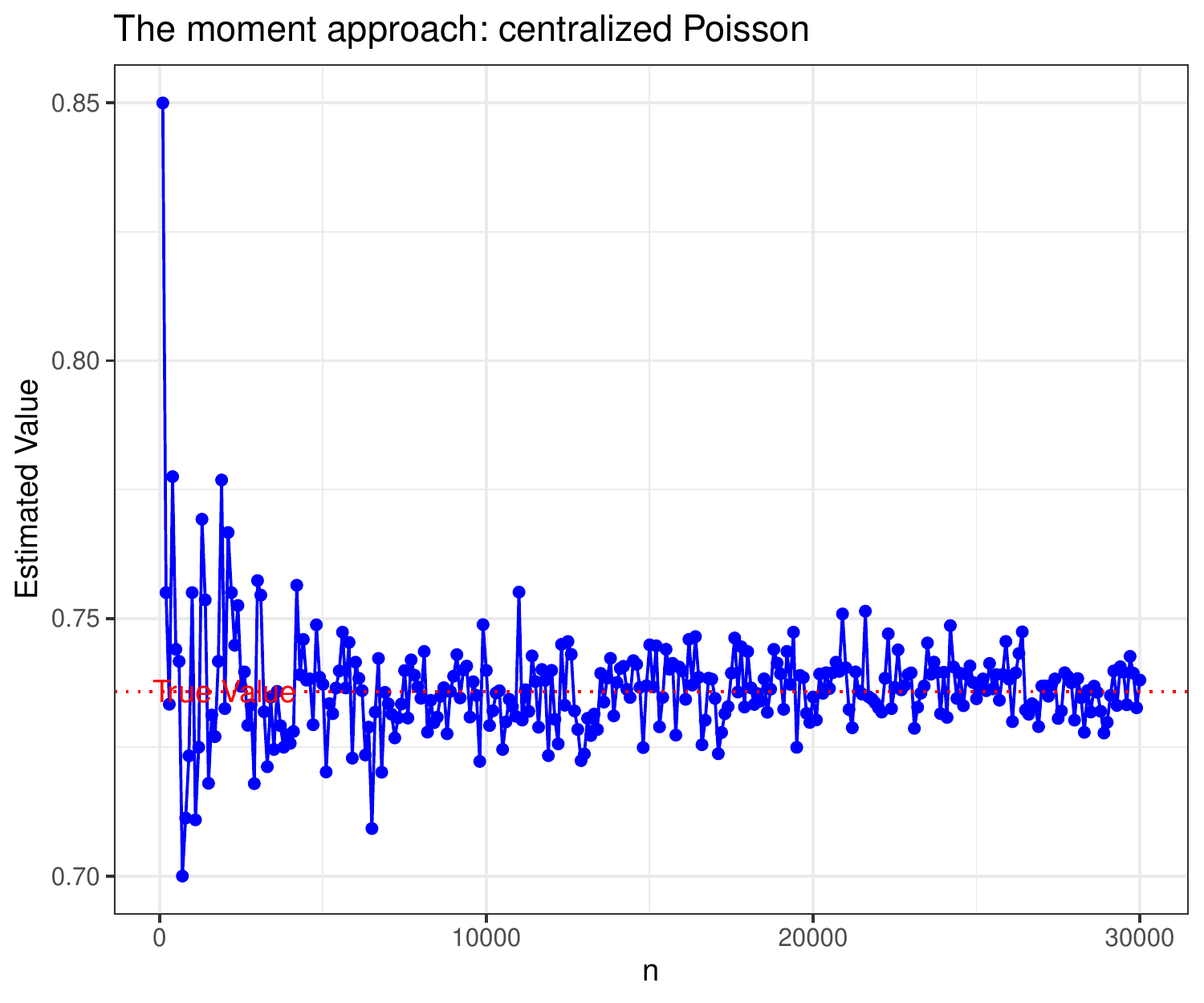}
        \caption{centralized Poisson}
        \label{figure5}
    \end{center}
    \endminipage
    \end{figure}

   For centralized negative binomial, and centralized Poisson ($\lambda = 1$) variable $X$, $\| X\|_{\varphi_1} =2.460938, 0.7357589,$ respectively. The Figure \ref{figure4} and Figure \ref{figure5} show the estimated value from different $n$ under estimate method (\ref{est3}) for the four distributions mentioned above.

The five figures mentioned above show litter bias between the estimated norm and true norm. It is worthy to note that the norm estimator for centralized negative binomial case has a peak point. This is caused by sub-exponential distributions having relatively heavy tails, and hence the norm estimation may not robust as that in sub-Gaussian under relatively small sample sizes.

Moreover, sub-Gaussian and sub-exponential parameter is extensible for random vectors with values in a normed space $(\mathcal{X},\|\cdot\|)$, we define \emph{norm-sub-Gaussian parameter} and \emph{norm-sub-exponential parameter}: The norm-sub-Gaussian parameter:
\begin{center}
${\left\| \boldsymbol{X} \right\|_{{\varphi _2}}} ={\sup }_{k \ge 1} {(k!)^{ - 1/(2k)}}{\left( {{\rm{E}}\|\boldsymbol{X}\|^{2k}} \right)^{1/(2k)}}$;
\end{center}
the norm-sub-exponential parameter:
\begin{center}
${\left\| \boldsymbol{X} \right\|_{{\varphi _1}}} ={\sup }_{k \ge 1} {(k!)^{ - 1/k}}{\left( {{\rm{E}}\|\boldsymbol{X}\|^{k}} \right)^{1/k}}.$
\end{center}
We denote $\boldsymbol{X} \sim\mathrm{nsubG}_{\varphi_1}(\sigma^2)$ and $\boldsymbol{X} \sim\mathrm{nsubG}_{\varphi_2}(\sigma^2)$ for $\sigma^2={\left\| \boldsymbol{X} \right\|_{{\varphi _2}}}~\text{and}~{\left\| \boldsymbol{X} \right\|_{{\varphi _1}}}$, respectively.

\section{Statistical Applications of Sub-Weibull Concentrations}

\subsection{Negative binomial regressions with heavy-tail covariates}
In statistical regression analysis, the responses $\{Y_{i}\}_{i=1}^n$ in linear regressions are assume to be continuous Gaussian variables. However, the category in classification or grouping may be infinite with index by the non-negative integers. The categorical variables is treated as countable responses for distinction categories or groups; sometimes it can be infinite. In practice, random count responses include the number of patients, the bacterium in the unit region, or stars in the sky and so on. The responses $\{Y_{i}\}_{i=1}^n$ with covariates $\{X_{i}\}_{i=1}^n$ belongs to generalized linear regressions. We consider i.i.d. random variables $\{(X_{i}, Y_{i})\}_{i=1}^n \sim (X,Y) \in \mathbb{R}^{p} \times \mathbb{N}$. By the methods of the maximum likelihood or the M-estimation, the estimator $\hat{\beta}_{n}$ is given by
\begin{equation}\label{eq:M-estimator}
\hat{\beta}_{n}:=\underset{\beta \in \mathbb{R}^{p}}{\arg \min } \frac{1}{n} \sum_{i=1}^{n} \ell (X_{i}^{\top}\beta, Y_{i}),
\end{equation}
where the loss function $\ell(\cdot, \cdot)$ is convex and twice differentiable in the first argument.

In  high-dimensional regressions, the dimension $\bm\beta$ may be growing with sample size $n$. When $\{Y_{i}\}_{i=1}^n$ belongs to the exponential family, \cite{Portnoy88} studied the asymptotic behavior of $\hat{\beta}_{n}$ in the generalized linear models (GLMs) as $p_n:=\mathrm{dim}(X)$ is increasing. In our study, we focus on the case that the covariates is $\operatorname{subW} (\theta)$ heavy-tailed for $\theta<1$.

The target vector $\beta^*:=\underset{\beta\in \mathbb{R}^{p}}{\arg \min }  \mathrm{E}\ell\left( X^{T}\beta, Y\right)$ is assumed to be the loss under the population expectation, comparing to \eqref{eq:M-estimator}.  Let $\dot\ell(u, y):=\left.\frac{\partial}{\partial t} \ell(t, y)\right|_{t=u}$, $\ddot\ell(u, y):=\left.\frac{\partial}{\partial t} \dot\ell(t, y)\right|_{t=u}$ and $C(u, y):=\sup _{|s-t| \leq u} \frac{\ddot\ell(s, y)}{\ddot\ell(t, y)}$. Finally, define the score function and Hessian matrix of the empirical loss function are $\hat{\mathcal{Z}}_{n}(\beta):=\frac{1}{n} \sum_{i=1}^{n} \dot\ell(X_{i}^{T} \beta, Y_{i})  X_{i}$ and $\hat{\mathcal{Q}}_{n}(\beta):=\frac{1}{n} \sum_{i=1}^{n} \ddot\ell(X_{i}^{T} \beta, Y_{i}) X_{i} X_{i}^{T}$, respectively. The population version of Hessian matrix is $\mathcal{Q}(\beta):=\mathrm{E}[\ddot\ell(X^{T}\beta, Y)  X X^{T}]$. The following so-called determining inequalities guarantee the $\ell_2$-error for the estimator obtained from the smooth M-estimator defined as \eqref{eq:M-estimator}.

\begin{lemma}[Corollary 3.1 in \cite{Kuchibhotla18D}] \label{norm on net general}
Let $\delta_{n}(\beta):=\frac{3}{2}\|[\hat{\mathcal{Q}}_{n}(\beta)]^{-1} \hat{\mathcal{Z}}_{n}(\beta)\|_{2}$ for ${\beta\in \mathbb{R}^{p}}$. If $\ell(\cdot, \cdot)$ is a twice differentiable function that is convex in the first
argument and for some $\beta^* \in \mathbb{R}^{p}$: $\max _{1 \leq i \leq n} C\left(\left\|X_{i}\right\|_{2} \delta_{n}(\beta^*), Y_{i}\right) \leq \frac{4}{3}$. Then there exists a vector $\hat{\beta}_{n} \in \mathbb{R}^{p}$ satisfying $\hat{\mathcal{Z}}_{n}(\hat{\beta}_{n})=0$ as the estimating equation of \eqref{eq:M-estimator},
$$
\frac{1}{2} \delta_{n}(\beta^*) \leq\|\hat{\beta}_{n}-\beta^*\|_{2} \leq \delta_{n}\left(\beta^*\right).
$$
\end{lemma}
Applications of Lemma \ref{norm on net general} in regression analysis is of special interest when $X$ is heavy tailed, i.e. the sub-Weibull index $\theta <1$. For the negative binomial regression (NBR) with the known dispersion parameter $k>0$, the loss
function is
\begin{equation}\label{eq:nbloss}
\ell(u, y)=-y u+(y+k) \log (k+e^{u}).
\end{equation}
 Thus we have $\dot\ell(u, y)=-{\frac{{k({y} - {e^u})}}{{k + {e^u}}}},~\ddot\ell(u, y)=\frac{k(y+k)e^{u}}{(k+e^{u})^{2}}$, see \cite{Zhang17} for details.

Further computation gives $C(u, y)=\sup _{|s-t| \leq u} \frac{ e^{s}(k+e^{t})^{2}}{(k+e^{s})^{2} e^{t}}$ and it implies that $C(u, y) \leq e^{3 u} .$ Therefore,
condition $\max _{1 \leq i \leq n} C\left(\left\|X_{i}\right\|_{2} \delta_{n}(\beta^*), Y_{i}\right) \leq \frac{4}{3}$ in Lemma \ref{norm on net general} leads to
\begin{center}
$
\max _{1 \leq i \leq n}\left\|X_{i}\right\|_{2} \delta_{n}\left(\beta^*\right) \leq \frac{\log (4 / 3)}{3}.
$
\end{center}
This condition need the assumption of the design space for  $\max _{1 \leq i \leq n}\left\|X_{i}\right\|_{2}$.

In NBR  with loss \eqref{eq:nbloss}, one has
\begin{center}
${\widehat {\cal Q}_n}(\beta^* ): = \frac{1}{n}\sum\limits_{i = 1}^n {{\frac{{({Y_i} + k)k{e^{X_i^{\top}\beta^* }}{X_i}X_i^{\top}}}{{{(k + {e^{X_i^{\top}\beta^* }})^2}}} }}$
and ${\widehat {\cal Z}_n}(\beta ^*): = \frac{-1}{n}\sum\limits_{i = 1}^n {\frac{{k({Y_i} - {e^{X_i^{\top}\beta ^* }}){X_i}}}{{k + {e^{X_i^{\top}\beta^* }}}}} $.
\end{center}
To guarantee that $\hat{\beta}_{n}$ approximates $\beta^*$ well, some regularity conditions are required.
\begin{itemize}
\item [\textbullet] (C.1):
For ${M_Y}, M_X>0$, assume $\mathop {\max }\limits_{1 \le i \le n} {\| {{Y_i} } \|_{{\psi _{1}}}} \le {M_Y}$ and the heavy-tailed covariates $\{X_{ik}\}$ are uniformly sub-Weibull with $\max\limits_{1 \le i \le n,1 \le k \le p}{\left\| {{X_{ik}}} \right\|_{{\psi _\theta }}} \le {M_X}$ for $0<\theta <1$.
\item [\textbullet] (C.2): The vector ${X_i}$ is sparse or bounded. Let ${\mathcal{F}_Y}: = \{ \mathop {\max }\limits_{1 \le i \le n} {\rm{E}}{Y_i} = \mathop {\max }\limits_{1 \le i \le n} {e^ {X_i^{\top}{ \beta ^*}}} \le B, \mathop {\max }\limits_{1 \le i \le n} \|{\bm X_i}\|_2 \le I_n\} $ with a slowly increasing function $I_n$, we have $\pr \{{\mathcal{F}_Y^c} \}  = {\varepsilon _n} \to 0$.
\end{itemize}
In addition, to bound $\mathop {\max }\limits_{1 \le i \le n,1 \le i \le k} |X_{ik}^{}| $, the sub-Weibull concentration determines:
\[\pr \big(\mathop {\max }\limits_{1 \le i \le n,1 \le i \le k} |X_{ik}| > t \big) \le np \pr (|X_{11}| > t) \le 2np{e^{ - {{(t/{{\left\| {{X_{11}}} \right\|}_{{\psi _\theta }}})}^\theta }}} \leq \delta  \Rightarrow t = {M_X}{\log ^{1/\theta }}(\frac{{2np}}{\delta }),\]
by using Corollary \ref{prop:Psub-W}. Hence, we define the event for the maximum designs:
\begin{equation*}
    {\mathcal{F}_{\max}} = \Big\{ {\mathop {\max }\limits_{1 \le i \le n,1 \le k \le p} |{X_{ik}}| \le {M_X}{\log ^{1/\theta }}(\frac{{2np}}{\delta })}\Big\} \cap {\mathcal{F}_Y}.
\end{equation*}
To make sure that the optimization in \eqref{eq:M-estimator} has a unique solution, we also require the minimal eigenvalue condition.
\begin{itemize}
\item [\textbullet] (C.3):
Suppose that $ {b^{\top}}{\rm{E}}(\hat{\mathcal{Q}}_{n}(\beta))b\ge C_{\min}$ is satisfied for all $b \in S^{p-1}$.
\end{itemize}
In the proof, to ensure that the random Hessian function has a non-singular eigenvalue, we define the event
 $${\mathcal{F}_1} = \left\{ \mathop {\max }\limits_{k,j}{\left| {\frac{1}{n}\sum\limits_{i = 1}^n {\left[ {\frac{{{Y_i}k{e^{X_i^{\top}\beta^* }}{X_{ik}}{X_{ij}}}}{{{(k + {e^{X_i^{\top}\beta^* }})^2}}} - {\rm{E}}\left( {\frac{{{Y_i}k{e^{X_i^{\top}\beta^* }}{X_{ik}}{X_{ij}}}}{{{(k + {e^{X_i^{\top}\beta^* }})^2}}}} \right)} \right]} } \right|}\le \frac{ C_{\min}}{4}\right\}$$
 $${\mathcal{F}_2} = \left\{ \mathop {\max }\limits_{k,j}\left| {\frac{1}{n}\sum\limits_{i = 1}^n {\left[ {\frac{{k{e^{X_i^{\top}\beta^* }}{X_{ik}}{X_{ij}}}}{{{(k + {e^{X_i^{\top}\beta^* }})^2}}} - {\rm{E}}\left( {\frac{{k{e^{X_i^{\top}\beta^* }}{X_{ik}}{X_{ij}}}}{{{(k + {e^{X_i^{\top}\beta^* }})^2}}}} \right)} \right]} } \right|\le \frac{ C_{\min}}{4}\right\}.$$
\begin{theorem}[Upper bound for $\ell_2$-error]\label{eq:L2uppercox}
In the NBR with loss \eqref{eq:nbloss} and $(C.1-C.3)$, let
\begin{equation*}
    M_{BX}=M_X + \frac{B}{{\log 2}},~~R_n := \frac{6 M_{BX} M_X}{C_{\min}} \left[ \sqrt{ \frac{2p}{n} \log \Big( \frac{2p}{\delta}\Big)} + \frac{1}{n} \sqrt{p \log \Big( \frac{2p}{\delta}\Big) } \right] \log^{1 / \theta} \Big( \frac{2 n p}{\delta}\Big),
\end{equation*}
and $\bm b := (k / n) M_X^2(1, \ldots ,1)^{\top}\in\mathbb{R}^n$. Under the event ${\mathcal{F}_1}\cap{\mathcal{F}_2}\cap{\mathcal{F}_{\max}}$, for any $0 < \delta < 1$, if the sample size $n$ satisfies
\begin{equation}\label{eq:CC}
    R_n I_n \le \frac{{\log (4/3)}}{3},
\end{equation}
Let $c_n:=e^{ - \frac{1}{4} (\frac{{n{t^2}}}{{2 M_X^4{{\log }^{4/\theta }}(\frac{{2np}}{\delta })M_{BX}^2}} \wedge \frac{{nt}}{{M_X^2{{\log }^{2/\theta }}(\frac{{2np}}{\delta })M_{BX}}})} + e^{ - (\frac{{{t^{\theta/2} }}}{{{{[4eC(\theta/2){{\left\|\bm b \right\|}_2}{L_n}(\theta/2 ,\bm b)]}^{\theta/2} }}} \wedge \frac{{{t^2}}}{{16{e^2}{C^2}(\theta/2)\left\|\bm b \right\|_2^2}})}$ with $t=C_{\min}/{4}$, then
\begin{equation*}
    \pr ({\| {{{\hat \beta }_n} - {\beta ^*}}\|_2} \le R_n)\ge 1-2p^2c_n-\delta-{\varepsilon _n}.
\end{equation*}
\end{theorem}

A few comment is made on this theorem. First, in order to get ${\| {{{\hat \beta }_n} - {\beta ^*}}\|_2} \xrightarrow{p} 0$, we need $p=o(n)$ under sample size restriction \eqref{eq:CC} with ${I_n} = o({\log ^{-1/\theta }}(np) \cdot {[{n^{ - {\rm{1}}}}p\log p{\rm{]}}^{-1/2}})$. Second, note that the ${\varepsilon _n}$ in provability $1-2p^2c_n-\delta-{\varepsilon _n}$ depends on the models size and the fluctuation of the design by the event ${\mathcal{F}_{\max}}$.

\subsection{Non-asymptotic Bai-Yin's theorem}
 In statistical machine learning, exponential decay tail probability is crucial to evaluate the finite-sample performance. Unlike Bai-Yin's law with the fourth-moment condition that leads to polynomial decay tail probability, under sub-Weibull conditions of data, we provide a exponential decay tail probability on the extreme eigenvalues of a $n \times p$ random matrix.

Let $\mathbf{A} = \mathbf{A}_{n,p}$ be an $n \times p$ random matrix whose entries are independent copies of a r.v. with zero mean, unit variance, and finite fourth moment. Suppose that the dimensions $n$ and $p$ both grow to infinity while the aspect ratio $p/n$ converges to a constant in $[0,1]$. Then Bai-Yin's law \citep{Bai1993} asserted that the standardized extreme eigenvalues satisfying
\begin{center}
$
  \frac{1}{{\sqrt n }}\lambda_{min}(\mathbf{A}) =1 - \sqrt {\frac{p}{n}}  + o\left( \sqrt {\frac{p}{n}} \right), \quad
 \frac{1}{{\sqrt n }} \lambda_{max}(\mathbf{A}) =1 + \sqrt {\frac{p}{n}}  + o\left(\sqrt {\frac{p}{n}} \right)\quad
  \text{a.s.}.
$
\end{center}
Next we introduce a special \emph{counting measure} for measuring the complexity of a certain set in some space. The $\mathcal{N}_\varepsilon$ is called an \emph{$\varepsilon$-net} of $K$ in $\mathbb{R}^n$ if
  $K$ can be covered by balls with centers in $K$ and radii $\varepsilon$ (under Euclidean distance). The \emph{covering number} $\mathcal{N}(K,\varepsilon )$ is defined by the smallest number of
  closed balls with centers in $K$ and radii $\varepsilon$ whose union covers $K$.

For purposes of studying random matrices, we need to extend the definition of sub-Weibull r.v. to sub-Weibull random vectors. The $n$-dimensional unit Euclidean sphere $S^{n-1}$, is denoted by ${S^{n - 1}} = \{ \boldsymbol{x} \in {\mathbb{R}^n}:{\left\| \boldsymbol{x} \right\|_2} = 1\}.$ We say that a random vector $\boldsymbol{X}$ in $\mathbb{R}^n$ is sub-Weibull
  if the one-dimensional marginals $\langle {\boldsymbol{X},\boldsymbol{a}} \rangle $ are sub-Weibull r.v.s for all $\boldsymbol{a} \in \mathbb{R}^n$.
The sub-Weibull norm of a random vector $\boldsymbol{X}$ is defined as
$\|\boldsymbol{X}\|_{\psi _\theta} := \sup_{\boldsymbol{a} \in S^{n-1}} \|\langle\boldsymbol{X},\boldsymbol{a} \rangle\|_{\psi _\theta}.$ Similarly, define the spectral norm for any $p \times p$ matrix $\mathbf{B}$ as $\big\|\mathbf{B}\big\|= \max_{||\textbf{x}||_2 =1}\big\|\mathbf{B}\boldsymbol{x}\big\|_2 =\sup_{\boldsymbol{x} \in S^{p-1}} |\langle {\mathbf{B}\boldsymbol{x},\boldsymbol{x}} \rangle|$. Spectral norm has many good properties, see \cite{Vershynin18} for details.

Furthermore, for simplicity, we assume that the rows in random matrices are isotropic random vectors. A random vector $\bm{Y}$ in $\mathbb{R}^{n}$ is called isotropic if $\mathrm{Var}(\bm{Y})=\mathbf{I}_p.$ Equivalently, $\bm{Y}$ is isotropic if
$\mathrm{E}[\langle \bm{Y}, \boldsymbol{a}\rangle^{2}]=\|\boldsymbol{a}\|_{2}^{2}~\text { for all } \boldsymbol{a} \in \mathbb{R}^{n}.$ In the non-asymptotic regime, Theorem 4.6.1 in \cite{Vershynin18} study the  upper and lower bounds of maximum (minimum) eigenvalues of random matrices with independent sub-Gaussian entries which are sampled from high-dimensional distributions. As an extension of Theorem 4.6.1 in \cite{Vershynin18}, the following result is a non-asymptotic versions of Bai-Yin's law for sub-Weibull entries, which is useful to estimate covariance matrices from heavy-tailed data [$\operatorname{subW} (\theta)$, $\theta<1$].

\begin{theorem}[{\rm{Non-asymptotic Bai-Yin's law}}] \label{NA-BY}
  Let $\mathbf{A}$ be an $n \times p$ matrix whose rows $\boldsymbol{A}_i$ are independent isotropic
  sub-Weibull random vectors in $\mathbb{R}^p$ with covariance matrix $\mathbf{I}_p$ and ${\max }_{1 \le i \le n} \|\bm{A}_i\|_{\psi _\theta}\le K$.
  Then for every $s \ge 0$, we have
\begin{equation*}
    {\pr} \bigg\{ \big\| \frac{1}{n}\mathbf{A}^{\top}\mathbf{A}-\mathbf{I}_p \big\|\leq H(cp + s^2, n ; \theta)\bigg\} \geq 1-2 e^{-s^2},
\end{equation*}
where
$$ H(t, n; \theta) := 2eKC(\theta/2 ){K_{\theta/2}[1+([({e\theta/2} )^{\theta/2}] \log 2)^{-\theta/2}}) \left[ \sqrt{\frac{t}{n}} + \left\{ \begin{array}{ll}
            A(\theta / 2) \frac{(\gamma^2 t)^{2 / \theta}}{n},  & \theta \leq 2\\
            B(\theta / 2) \frac{(\gamma^2 t)^{2 / \theta}}{{n}^{1/\theta }}, &\theta>2
        \end{array} \right. \right],$$
where $K_{\alpha}:=2^{1 / \alpha}$ if $\alpha \in(0,1)$ and $K_{\alpha}=1$ if $\alpha \geq 1$; $A(\theta /2)$, $B(\theta /2)$ and $C(\theta /2)$ defined in Theorem \ref{thm:SumNewOrliczVex}(a).

Moreover, the concentration inequality for extreme eigenvalues hold for $c \geq {n \log 9/p}$
\begin{equation}\label{eq:P2}
     \pr\Big\{\sqrt{1-H^2(cp + s^2, n ; \theta)} \le \frac{\lambda_{min}(\mathbf{A})}{\sqrt{n}} \le \frac{\lambda_{max}(\mathbf{A})}{\sqrt{n}} \le \sqrt{1+H^2(cp + s^2, n ; \theta)}\Big\} \ge 1-2 e^{-s^2}.
\end{equation}
\end{theorem}

\subsection{General Log-truncated Z-estimators and sub-Weibull type robust estimators}
Motivated from log-truncated loss in \cite{chen2021generalized,Xu22}, we study the almost surely continuous and non-decreasing function $\varphi^c: \mathbb{R} \rightarrow \mathbb{R}$ for truncating the original score function
\begin{align}\label{eq:CANTON}
	-\log \left[1-x+c(|x|)\right]\leq \varphi^c(x) \leq \log \left[1+x+c(|x|)\right], \quad \forall x \in \mathbb{R}
\end{align}
where $c(|x|)>0$ is a \emph{high-order function} \citep{Xu22} of $|x|$ which is to be specified. For example, a plausible choose for $\varphi^c(x)$ in \eqref{eq:CANTON} should have following form
	\begin{align}\label{eq:CANTONsp}
		{\varphi ^c}(x)&=\log \left[ {1 + x + c(|x|)} \right]{\rm{1(}}x \ge 0{\rm{)}} - \log \left[ {1 - x + c(|x|)} \right]{\rm{1(}}x \le 0{\rm{)}}\nonumber\\
		&={\rm{sign}}(x)\log (1 + |x| +c(|x|)).
	\end{align}
	For \eqref{eq:CANTONsp}, we get ${\varphi ^c}(x) \approx x$ for sufficiently smaller $x$ and ${\varphi ^c}(x) \ll  x$ for larger $x$. Under \eqref{eq:CANTON}, now we show that $c(|x|)$ must obey a key inequality. For all $x \in \mathbb{R}$, it suffices to verify
$ - \log [1 - x + {\rm{c}}(|x|)] \le \log [1 + x + {\rm{c}}(|x|)]$,
which is equivalent to check 
   $\log\left[ {\left( {1 + {\rm{c}}(|x|) + x} \right)\left( {1 + {\rm{c}}(|x|) - x} \right)} \right] \ge 0$, namely
${\left( {1 + {\rm{c}}(|x|)} \right)^2} - {x^2} \ge 1 \Leftrightarrow {\rm{c}}(|x|) \ge \sqrt {1 + {x^2}}  - 1.$ 

For independent r.v.s $\{X_{i}\}_{i=1}^n$, using the score function \eqref{eq:CANTONsp}, we define the score function of data
\begin{center}
${{\hat Z}_{\alpha_n} }(\theta ) = \frac{1}{{n\alpha_n }}\sum_{i = 1}^n \varphi^c \left[ {\alpha_n \left( {{X_i} - {\theta }} \right)} \right]$ for any $\theta \in \mathbb{R}$. 
\end{center}
Then the influence of the heavy-tailed outliers is weaken by $\varphi^c \left[ {\alpha_n \left( {{X_i} - {\theta }} \right)} \right]$ by choosing an optimal $\alpha_n$.
	We aim to estimate the \emph{average mean}: $\mu_n: = \frac{1}{n}\sum_{i = 1}^n {\rm{E}}{X_i}$ for non-i.i.d. samples $\{X_i\}_{i=1}^{n}$. Define the $Z$-estimator ${{\hat \theta }_{\alpha_n} }$ as
	\begin{align}\label{eq:CANTONZ}
		{{\hat \theta }_{\alpha_n} } \in \{ \theta \in \mathbb{R} :{{\hat Z}_{\alpha_n} }(\theta ) = 0\} ,
	\end{align}
	where $\alpha_n$ is the tuning parameter (will be determined later).

To guarantee consistency for log-truncated Z-estimators \eqref{eq:CANTONZ}, we require following assumptions of $c(\cdot)$.
	\begin{itemize}
		\item [\textbullet] (C.1): For a constant $c_2>1$, the $c(x)$ satisfies \emph{weak triangle inequality} and \emph{scaling property},
		$${\rm(C.1.1)}: c(x+y)\le c_2[c(x)+c(y)],~~{\rm(C.1.2)}:c(tx) \le f(t)c(x)$$
		for $f(t)$ satisfies (C.1.3): $f(t)$ and $f(t)/|t|$ are non-constant increasing functions and $\mathop {\lim }\limits_{t \to 0} f(t)/|t| = 0$.
	\end{itemize}
	
\begin{remark}
	Note that  $|x| \ge \sqrt {1 + {x^2}}  - 1$ and we could put
	$c(|x|)=|x|$. However, $c(|x|)=|x|$ does not satisfy (C.1.3) since $f(t)=|t|$ and $f(t)/|t|$ are constant functions of $t$.
\end{remark}	

In the following theorem, we establish the finite sample confidence interval and the convergence rate of the estimator ${{\hat \theta }_{\alpha_n} }$.

\begin{theorem}\label{eq:cantoniG}
Let $\{X_i\}_{i=1}^{n}$ be  independent samples drawn from an unknown probability distribution $\{P_i\}_{i=1}^{n}$ on $\mathbb{R}$. Consider the estimator $\hat \theta_{\alpha_n}$ defined as \eqref{eq:CANTONZ} with (C.1), ${\alpha_n} \to 0$ and $\frac{1}{{n}}\sum_{i = 1}^n {{\rm{E}}[c({X_i} - \theta )]}=O(1)$. Let $B_n^ + (\theta ) = \mu_n  - \theta  + \frac{1}{{n\alpha_n }}\sum_{i = 1}^n {\rm{E}}[c\left( {\alpha_n ({X_i} - \theta )} \right)] {\rm{ + }}\frac{{\log ({\delta ^{ - 1}})}}{n\alpha_n}$ and $B_n^ - (\theta ) = \mu_n  - \theta  - \frac{1}{{n\alpha_n }}\sum_{i = 1}^n {{\rm{E}}[c\left( {\alpha_n ({X_i} - \theta )} \right)]}  - \frac{{\log ({\delta ^{ - 1}})}}{n\alpha_n }$. Let $\theta _{+}$ be the smallest solution of the equation $B_n^ + (\theta ) =0$ and $\theta_{-}$ be the largest solution of $B_n^ - (\theta )=0$.\\
{\rm{(a)}}. We have with the $(1 - 2\delta)$-confidence intervals
\[
\pr(B_n^ - (\theta ) < {{\hat Z}_{\alpha_n} }(\theta ) < B_n^ + (\theta )) \ge 1 - 2\delta, \qquad \pr(\theta_{-} \leq \hat{\theta}_{\alpha_n} \leq \theta_{+}) \ge 1 - 2\delta,
\]
for any $\delta \in(0,1 / 2)$ satisfies the sample condition:
\begin{align}\label{eq:sample condition}
	\frac{1}{{n\alpha_n }}\sum\limits_{i = 1}^n {{\rm{E}}[c\left( {\alpha_n {X_i} - \alpha_n [\mu_n  \pm  {d_n}(c)]} \right)]} +\frac{{\log ({\delta ^{ - 1}})}}{n\alpha_n } < {d_n}(c),
\end{align}
where ${d_n}(c)$  is a constant such that $B_n^ \pm (\mu_n\pm d_n(c) )<0$.\\
{\rm{(b)}}. Moreover, picking $\alpha_n  \ge {f^{ - 1}}\left( {\frac{{\log ({\delta ^{ - 1}})}}{{c_2\sum\nolimits_{i = 1}^n {\rm{E}}[c\left( {{X_i} - \mu_n } \right)]} }} \right)$, one has
		\begin{align}\label{eq:con}
\pr \left(| {\hat \theta_{\alpha_n}-\mu_n } | \le \left| {g_{\alpha_n} ^{ - 1}\left\{-{\frac{{2\log ({\delta ^{ - 1}})}}{n\alpha_n }}\right\}} \right|\right)\ge 1-2\delta,~\text{with}~{g_{\alpha_n} }(t) := t + \frac{{{c_2}}}{\alpha_n }c\left( {\alpha_n t} \right).
		\end{align}
	\end{theorem}
The \eqref{eq:con} in Theorem \ref{eq:cantoniG} is a fundamental extension of Lemma 2.1 (see Theorem 16 in \cite{lerasle2019lecture}) with $c(x)=x^2/2$ from i.i.d. sample to independent sample. Let $c(x)=|x|^{\beta}/{\beta}$ , for i.i.d. sample, Theorem \ref{eq:cantoniG} implies Lemmas 2.3, 2.4 and Theorem 2.1 in \cite{chen2021generalized}. The  $\alpha_n  \ge {f^{ - 1}}\left( {\frac{{\log ({\delta ^{ - 1}})}}{{c_2\sum\nolimits_{i = 1}^n {\rm{E}}[c\left( {{X_i} - \mu_n } \right)]} }} \right)$ in Theorem \ref{eq:cantoniG}(b) gives a theoretical guarantee for choosing the tuning parameter $\alpha_n$.
\begin{proposition}[Theorem 2.1 in \cite{chen2021generalized}]\label{coro:cp}
	Let $\{X_{i}\}_{i=1}^{n}$ be a sequence of i.i.d. samples drawn from an unknown probability distribution on $\mathbb{R}$. We assume ${\rm{E}}\left|X_{1}\right|^{\beta}<\infty$ for a certain ${\beta} \in(1,2]$ and denote $\mu={\rm{E}}\left[X_{1}\right],~v_{\beta}={\rm{E}}\left|X_{1}-\mu\right|^{\beta}$. Given any $\epsilon \in(0,1 / 2)$ and positive integer $n \geq\left(\frac{2 v_{\beta}+1}{\beta}\right)^{\frac{\beta}{{\beta}-1}} \frac{2 {\beta} \log \left(\epsilon^{-1}\right)}{v_{\beta}}$, let ${\alpha_n}=\frac{1}{2}(\frac{2 {\beta} \log (\epsilon^{-1})}{n v_{\beta}})^{\frac{1}{\beta}} .$ Then, with probability at least $1-2\epsilon$,
\begin{align}\label{eq:chen}
	| {\hat \theta_{\alpha_n}-\mu } |\le {2\left(\frac{2 {\beta} \log (\epsilon^{-1})}{n}\right)^{\frac{{\beta}-1}{\beta}} v_{\beta}^{\frac{1}{\beta}}}\left[{{\beta}-\left(\frac{2 {\beta} \log (\epsilon^{-1})}{n v_{\beta}}\right)^{\frac{{\beta}-1}{\beta}}}\right]^{-1}=O\left(n^{-\frac{{\beta}-1}{\beta}}\right).
\end{align}
\end{proposition}

Comparing to the convergence rate in \eqref{eq:chen}, put $O(n^{-\frac{{\beta}-1}{\beta}})=O({n^{ - 1/\theta }})$ for $\theta  > 2$. It implies
\[\beta^{-1 } + \theta^{-1 }  = 1,~(\theta  \ge 2~\text{or}~0 < \beta  \le 2).\]

For example, let us deal with the Pareto distribution ${\rm{Pareto}}({\alpha},k)$ with shape parameter $\alpha>0$ and scale parameter $k>0$, and the density function is $f(x)=\frac{\alpha k^{\alpha}}{x^{\alpha+1}}\cdot {\rm1}_{\{x \in [k, \infty )\}}$. For ${\alpha} \le 2$, ${\rm{Pareto}}({\alpha},k)$ has infinite variance, and it does not belong to the sub-Weibull distribution, so do the sample mean of i.i.d. Pareto distributed data. Proposition \ref{coro:cp} shows that the estimator error for robust mean estimator enjoys sub-Weibull concentration as presented in Proposition \ref{thm:orlicz_concentra}, without finite sub-Weibull norm assumption of data. With the Weibull-tailed behavior, it motivates us to define general sub-Weibull estimators having the non-parametric convergence rate $O({n^{ - 1/\theta }})$ in Proposition \ref{thm:orlicz_concentra} for $\theta  > 2$, even if the data do  not have finite sub-Weibull norm.
\begin{definition}[Sub-Weibull estimators]\label{def:estimators}
An estimator $\hat{\mu}:=\hat{\mu}(X_1,\cdots,X_n)$ based on i.i.d. samples $\{X_i\}_{i=1}^n$ from an unknown probability distribution ${P}$ with mean $\mu_{P}$, is called $(A, B, C)$-$\operatorname{subW}(\theta)$ if
$$
\forall t \in(0, A), \quad {\pr}(|\hat{\mu}-\mu_{P}|\le B{( {{t}/{n}} )^{1/\theta }}) \ge 1-C e^{-t}.
$$
\end{definition}
For example, in Proposition \ref{coro:cp}, ${\hat \theta_{\alpha_n}}$ is $(\infty, B, 1)$-$\operatorname{subW}({\frac{{\beta}}{\beta-1}})$ with  $B\sim {2\left({2 {\beta} \log (\epsilon^{-1})}\right)^{\frac{{\beta}-1}{\beta}} v_{\beta}^{\frac{1}{\beta}}}$ in Definition \ref{def:estimators}.
When $\theta=2$, \cite{devroye2016sub} defined sub-Gaussian estimators (includes Median of means and Catoni's estimators) for certain heavy-tailed distributions and discussed the nonexistence of sub-Gaussian mean estimators under $\beta$-moment condition for the data ($\beta \in (1,2)$).

\section{Conclusions}
Concentration inequalities are far-reaching useful in high-dimensional statistical inferences and machine learnings. They can facilitate various explicit non-asymptotic confidence intervals as a function of the sample size and model dimension.  

Future research includes sharper version of Theorem \ref{thm:sub-con} that is crucial to construct non-asymptotic and data-driven confidence intervals for the sub-Weibull sample mean. Although we have obtained sharper upper bounds for sub-Weibull concentrations, the lower bounds on tail probabilities are also important in some statistical applications \citep{Zhangan2020}. Developing non-asymptotic and sharp lower tail bounds of Weibull r.v.s is left for further study. For negative binomial concentration inequalities in Corollary \ref{remark2}, it is of interesting to study concentration inequalities of COM-negative binomial distributions (see \cite{zhang2018com}).

\section{Acknowledgement}

This work is supported n part by National Natural Science Foundation of China Grant (12101630) and  the University of Macau under UM Macao Talent Programme (UMMTP-2020-01). This work is also supported in part by the Key Project of Natural Science Foundation of Anhui Province Colleges and Universities (KJ2021A1034), Key Scientific Research Project of Chaohu University (XLZ-202105). The authors thank Guang Cheng for the discussion about doing statistical inference in the non-asymptotic way and Arun Kumar Kuchibhotla for his help about the proof of Theorem \ref{thm:SumNewOrliczVex}. The authors also thank Xiaowei
Yang for his helpful comments on Theorem \ref{eq:cantoniG}.


\newpage
\section{Appendix}

\begin{proof}[\hl{Proof of Corollary} 
 \ref{gote}]
$\phi_{|X|^{\theta}}(t)$ is continuous for $t$ a neighborhood of zero, by the definition,
$ 2 \ge \E e^{ (|X| / \| X\|_{\psi_{\theta}} )^{\theta}} = m_{|X|^{\theta}} \big( \| X\|_{\psi_{\theta}}^{-\theta}\big).$
    Since $|X|^{\theta} > 0$, the MGF $m_{|X|^{\theta}}(t)$ is monotonic increasing. Hence, inverse function $m_{|X|^{\theta}}^{-1}(t)$ exists and $\| X\|_{\psi_{\theta}}^{- \theta} = m_{|X|^{\theta}}^{-1}(2).$
    So $\| X\|_{\psi_{\theta}} = \big( m_{|X|^{\theta}}^{-1}(2) \big)^{- 1 / \theta}$.
\end{proof}

\subsection{}

\begin{proof}[\hl{Proof of Corollary} \ref{remark2}]
    The first inequality is the direct application of \eqref{exp} by observing that for any constant $a \in \mathbb{R}$, and r.v. $Y$ with $\| Y\|_{\psi_1} < \infty$, $ \| a Y\|_{\psi_1} = |a| \| Y \|_{\psi_1}$, $\|Y +a \|_{\psi_1} \leq \|Y\|_{\psi_1} + \|a\|_{\psi_1} = \|Y\|_{\psi_1} + |a| / \log 2$ and $\|X +a \|_{\psi_1}^2 \leq (\|X\|_{\psi_1} + |a| / \log 2)^2$. The second inequality is obtained from \eqref{exp} by considering two rate in $(\frac{t^{2}}{\sum_{i=1}^{n} 2\|Y_{i}\|_{\psi_{1}}^{2}} \wedge \frac{t}{\max _{1 \leq i \leq n}\|Y_{i}\|_{\psi_{1}}})$ separately. For \eqref{eq:NBRcon}, we only need to note that
    \begin{equation*}
        \begin{aligned}
            \scalebox{0.90}{$\| {Y_i} \|_{\psi_1}  = \inf \{t>0: \E e^{{Y_i} / t} \leq 2\} =\inf\{t>0: \left( \frac{1 - q_i}{1 - q_i e^{1 / t}}\right)^{k_i} \leq 2\}  = \left[ \log \frac{1 - (1 - q_i) / \sqrt[k_i]{2}}{q_i}\right]^{-1}$}.
        \end{aligned}
    \end{equation*}
    
    Then the third inequality is obtained by the first inequality and the definition of $a (\mu_i, k_i)$.
\end{proof}

\subsection{}

\begin{proof}[\hl{Proof of Corollary} \ref{prop:Psub-W}]
    The first and second part of this proposition were shown in Lemma 2.1 of \cite{Zajkowski19}. For the third result, using the bounds of Gamma function [see \cite{jameson2015a}]:
\begin{center}
$\sqrt {2\pi } {x^{x - (1/2)}}{e^{ - x}} \le \Gamma (x) \le [\sqrt {2\pi } {x^{x - (1/2)}}{e^{ - x}}] \cdot {e^{1/(12x)}},(x > 0),$
\end{center}
 it gives
\begin{align*}
{({\rm{E}}|X{|^k})^{1/k}}&  \scalebox{0.99}{$\le {\left\{ {2{\|X\|_{\varphi_{\theta}}^k}\left(  \frac{k}{\theta } \right)[\sqrt {2\pi } {{\left( {k/\theta } \right)}^{\frac{k}{\theta } - \frac{1}{2}}}{e^{ - \frac{{11k}}{{12\theta }}}}]} \right\}^{1/k}} = {( {\frac{{2\sqrt {2\pi } }}{\theta }})^{1/k}}{\{ {{{\left( \frac{k}{\theta } \right)}^{\frac{k}{\theta } + \frac{1}{2}}}{e^{ - \frac{{11k}}{{12\theta }}}}} \}^{1/k}}\|X\|_{\varphi_{\theta}}$}\\
& \scalebox{0.99}{$= {( {\frac{{2\sqrt {2\pi } }}{\theta }} )^{1/k}}{\left( {k/\theta } \right)^{\frac{1}{\theta } + \frac{1}{{2k}}}}{e^{ - \frac{{11}}{{12\theta }}}}\|X\|_{\varphi_{\theta}} \le {C_\theta }{( {\theta {e^{11/12}}})^{ - 1/\theta }}\|X\|_{\varphi_{\theta}}{k^{1/\theta }}$}.
\end{align*}
\end{proof}

\subsection{}

\begin{proof}[\hl{Proof of Corollary} \ref{lem:2theta}]
    By the definition of $\psi_{\theta}$-norm, ${\rm{E}}\exp \{ |{X}/{\| {{X}}\|_{{\psi _{\theta}}}}|^{\theta}\}\le 2$. Then \linebreak${\rm{E}}\exp \{ ||X|^r/{\| {{X}}\|^r_{{\psi _{\theta}}}}|^{\theta/r}\}\le 2.$ The result $|X|^r \sim \operatorname{subW}(\theta/r)$ follows by the definition of $\psi_{\theta}$-norm again. Moreover,
    \begin{align*}
        \|X\|_{\psi_{\theta}}:&=\inf \{C\in(0, \infty): ~ \mathrm E[\exp(|X|^{\theta}/C^{\theta})]\leq 2\}\\
&=[\inf \{C^r\in(0, \infty): ~ \mathrm E[\exp\{ ||X|^r/C^r|^{\theta/r}\}]\leq 2\}]^{1/r}=\left\||X|^{r}\right\|^{1/r}_{\psi_{\theta/r}},
    \end{align*}
which verifies \eqref{eq:2theta}. If  $X \sim \operatorname{subW}(r\theta)$, then ${\rm{E}}\exp \{ |X^r/{\| {{X}}\|^r_{{\psi _{r\theta}}}}|^{\theta}\}={\rm{E}}\exp \{ |X/{\| {{X}}\|_{{\psi _{r\theta}}}}|^{r\theta}\}\le 2$, which means that $X^r \sim \operatorname{subW}(\theta)$ with
    \begin{align*}
        \|X\|_{\psi_{r\theta}}:&=\inf \{C\in(0, \infty): ~ \mathrm E[\exp(|X|^{{r\theta}}/C^{{r\theta}})]\leq 2\}\\
&=[\inf \{C^r\in(0, \infty): ~ \mathrm E[\exp\{ ||X|^r/C^r|^{\theta}\}]\leq 2\}]^{1/r}=\left\||X|^{r}\right\|^{1/r}_{\psi_{\theta}}.
    \end{align*}
\end{proof}

\subsection{}

\begin{proof}[\hl{Proof of Corollary} \ref{lem2}]
    Set $\Delta := \sup_{p \ge 2} \frac{\| X\|_p}{\sqrt{p} + L p^{1/\theta}}$ so that $
        \|X\|_p \leq \Delta \sqrt{p} + L \Delta p^{1 / \theta}
   $ holds for all $p \geq 2$. By Markov's inequality for $t$-th moment $(t \ge 2)$, we have
    \begin{equation*}
        \begin{aligned}
             \pr \left( |X| \geq e \Delta \sqrt{t} + e L \Delta t^{1 / \theta} \right) & \le {\left( {\frac{{||X||_t}}{{e\Delta [\sqrt t  + L{t^{1/\theta }}]}}} \right)^t} \leq e^{-t},~[\text{By the definition of}~\Delta].
        \end{aligned}
    \end{equation*}
    
    So, for any $t \ge 2$,
    \begin{equation}\label{S.2.7}
        \pr \left( |X| \geq e \Delta \sqrt{t} + e L \Delta t^{1 / \theta} \right) \leq  {e^{ - t}}.
    \end{equation}
    
    Note the definition of $\Delta$ shows $
        \|X\|_t \leq \Delta \sqrt{t} + L \Delta t^{1 / \theta}
   $ holds for all $t \geq 2$ and assumption $\|X\|_t \leq C_1 \sqrt{t} + C_2 t^{1/\theta}$  for all $t \geq 2$. It gives $e \Delta \sqrt{t} + e L \Delta t^{1 / \theta} \leq  eC_1 \sqrt{t} + eC_2 t^{1/\theta}$. This inequality with (\ref{S.2.7}) gives
    \begin{equation}\label{eq:C1C2}
        \pr \left( |X| \geq e C_1 \sqrt{t} + eC_2 t^{1/\theta} \right) \leq 1\{ 0 < t < 2\}  + {e^{ - t}}\{ t \ge 2\} ,~~\forall t > 0.
    \end{equation}
    
    Take $K = k^{2 / \theta} C_2 / (k C_1)$, and define $\delta_k := ke C_1$ for a certain constant $k>1$,
    \begin{equation*}
        \begin{aligned}
             &~~~~\mathrm{E} \left[ \Psi_{\theta, K} \left( \frac{|X|}{ \delta_k }\right)\right] = \int_0^{\infty} \pr \left( |X| \geq \delta_k \Psi_{\theta, K}^{-1} (s)\right) ds \\
                         & = \int_0^{\infty} \pr ( |X| \geq ke C_1 \sqrt{\log (1 + s)} + ke C_1 K [\log (1 + s)]^{1/\theta}) ds \\
            & = \int_0^{\infty} \pr ( |X| \geq e C_1 \sqrt{\log (1 + s)^{k^2}} + e C_2  [\log (1 + s)^{k^2}]^{1/\theta}) ds \\
         [\text{By}~\eqref{eq:C1C2}]~   &  \le \int_{0 < {k^2}\log (1 + s) < 2} {ds}  + \int_{{k^2}\log (1 + s) \ge 2} {\exp \left\{ { - {k^2}\log (1 + s)} \right\}ds} \\
         &  \le \int_0^{{e^{{2k^{ - 2}}}} - 1} {dt}  + \int_{{e^{{2k^{ - 2}}}} - 1}^\infty  {\frac{{dt}}{{{{(1 + t)}^{{k^2}}}}}}\\
           &= {e^{{2k^{ - 2}}}} - 1 + \frac{{{{(1 + t)}^{1 - {k^2}}}}}{{1 - {k^2}}}\left. {} \right|_{{e^{{2k^{ - 2}}}} - 1}^\infty  = {e^{{2k^{ - 2}}}} - 1 + \frac{{{e^{2(1 - {k^2})/{k^2}}}}}{{{k^2} - 1}} \le 1.
        \end{aligned}
    \end{equation*}
    
    Therefore, $\|X\|_{\Psi_{\theta, K}} \leq  \gamma e C_1$ with $\gamma$ defined as the smallest solution of the inequality $\{ k > 1:{e^{{2k^{ - 2}}}} - 1 + \frac{{{e^{2(1 - {k^2})/{k^2}}}}}{{{k^2} - 1}} \le 1\}$. An approximate solution is $\gamma \approx 1.78$.
\end{proof}
 
\subsection{}\label{thm1.proof}
    The main idea in the proof is by the sharper estimates of the GBO norm of the sum of \emph{symmetric} r.v.s.
\begin{proof}[\hl{Proof of Theorem} \ref{thm:SumNewOrliczVex}]

\item (a) Without loss of generality, we assume $\|X_i\|_{\psi_\theta} = 1$. Define $Y_i := \big( |X_i| - ( \log 2)^{1 / \theta}\big)_+$, then it is easy to check that $\pr (|X_i| \geq t) \leq 2 e^{- t^{\theta}}$ implies $\pr (Y_i \geq t) \leq e^{- t^{\theta}}$. For independent Rademacher r.v. $\{\varepsilon_i\}_{i = 1}^n$, the symmetrization inequality gives $\left\| \sum_{i = 1}^n w_i X_i\right\|_p \leq 2 \left\| \sum_{i = 1}^n \varepsilon_i w_i X_i\right\|_p.$
Note that $\varepsilon_i X_i$ is identically distributed as $\varepsilon_i |X_i|$,
\begin{align}\label{eq:symmetrization}
\| \sum_{i = 1}^n w_i X_i\|_p & \leq 2 \| \sum_{i = 1}^n \varepsilon_i w_i |X_i|\|_p \leq 2 \| \sum_{i = 1}^n \varepsilon_i w_i  \big( Y_i + ( \log 2)^{1 / \theta}\big) \|_p \nonumber\\
 & \leq 2 \| \sum_{i = 1}^n \varepsilon_i w_i Y_i\|_p + 2 ( \log 2)^{1 / \theta} \| \sum_{i = 1}^n \varepsilon_i w_i\|_p \nonumber\\
[\text{Khinchin-Kahane inequality}]  & \leq 2 \| \sum_{i = 1}^n \varepsilon_i w_i Y_i\|_p + 2 ( \log 2)^{1 / \theta} \left(\frac{p-1}{2-1}\right)^{1 / 2}\| \sum_{i = 1}^n \varepsilon_i w_i\|_2 \nonumber\\
 & < 2\| \sum_{i = 1}^n \varepsilon_i w_i Y_i\|_p + 2 ( \log 2)^{1 / \theta} \sqrt{p}(\E(\sum_{i=1}^{n} \varepsilon_{i} w_{i})^{2})^{1 / 2}\nonumber\\
 [\left\{\varepsilon_{i}\right\}_{i=1}^{n} \text { are independent}]& =2 \| \sum_{i = 1}^n \varepsilon_i w_i Y_i\|_p + 2 ( \log 2)^{1 / \theta} \sqrt{p} \|\bm{w}\|_2.
\end{align}

From Lemma \ref{lem3}, we are going to handle the first term in \eqref{eq:symmetrization} with the sum of symmetric r.v.s. Since $\pr (Y_i \geq t) \leq e^{- t^{\theta}}$, then
\begin{center}
$ \| \sum_{i = 1}^n \varepsilon_i w_i Y_i\|_p = \| \sum_{i = 1}^n w_i Z_i \|_p, \qquad Z_i:=\varepsilon_i Y_i$
\end{center}
for symmetric independent r.v.s $\{Z_i\}_{i = 1}^n$ satisfying $|Z_i| \overset{d}{=}Y_i$ and $\pr (Z_i \geq t) = e^{- t^{\theta}}$ for all $t \geq 0$.

Next, we proceed the proof by checking the moment conditions in Corollary \ref{lem2}.

\emph{\hl{Case} 
 $\theta \leq 1$:} $N(t) = t^{\theta}$ is concave for $\theta \leq 1$. From Lemmas \ref{lem3} and \ref{lem4} (a), for $p \geq 2$,
\begin{equation*}
    \begin{aligned}
       &~~~~\Big\| \sum_{i = 1}^n w_i Z_i\Big\|_p  \leq e \inf \bigg\{ t > 0 \, : \, \sum_{i = 1}^n \log \phi_p \Big( e^{-2} \Big( \frac{w_i e^2}{t}\Big)Z_i\Big) \leq p\bigg\} \\
        & \leq e \inf \bigg\{ t > 0 \, : \, \sum_{i = 1}^n p M_{p, Z_i}  \Big( \frac{w_i e^2}{t}\Big) \leq p \bigg\} \\
        & = e \inf \bigg\{ t > 0 \, : \, \sum_{i = 1}^n \bigg[\Big\{\Big( \frac{w_i e^2 }{t}\Big)^p \|Z_i\|_p^p \Big\} \vee \Big\{p \Big( \frac{ w_i e^2}{t}\Big)^2 \|Z_i\|_2^2 \Big\}\bigg] \leq p \bigg\} \\
        & \leq e \inf \bigg\{ t > 0 : \Gamma \Big( \frac{p}{\theta} + 1\Big) \frac{e^{2p}}{t^p} \|\bm{w}\|_p^p \leq 1 \bigg\} + e \inf \bigg\{ t > 0 : p \Gamma \Big( \frac{2}{\theta} + 1\Big) \frac{e^{4}}{t^2} \|\bm{w}\|_2^2 \leq 1 \bigg\},
    \end{aligned}
\end{equation*}
where the last inequality we use $\| Z_i \|_p^p = \int_0^{\infty} p t^{p - 1} \pr (|Z_i| \geq t) \, dt \leq \int_0^{\infty} p t^{p - 1} e^{- t^{\theta}} \, dt = p \Gamma \left( \frac{p}{\theta} + 1\right).$
Hence
\begin{center}
    $\| \sum_{i = 1}^n w_i Z_i \|_p  \leq e^3 \left[  \Gamma^{1 / p} \left( \frac{p}{\theta} + 1\right)  \| \bm{w}\|_p+  \sqrt{p} \Gamma^{1 / 2} \left( \frac{2}{\theta} + 1\right) \| \bm{w}\|_2\right],$
\end{center}
and
\begin{equation*}
    \begin{aligned}
        \| \sum_{i =  1}^n w_i X_i\|_p  &\leq 2 e^3  \Big[  \Gamma^{1 / p} \Big( \frac{p}{\theta} + 1\Big)  \| \bm{w}\|_p + \sqrt{p} \Gamma^{1 / 2} \Big( \frac{2}{\theta} + 1\Big) \| \bm{w}\|_2 \Big] + 2 (\log 2)^{1 / \theta} \sqrt{p} \|\bm{w}\|_2 \\
        & = 2 e^3  \Gamma^{1 / p} \Big( \frac{p}{\theta} + 1\Big)  \| \bm{w}\|_p +  2\Big[(\log 2)^{1 / \theta} + e^3 \Gamma^{1 / 2} \Big( \frac{2}{\theta} + 1\Big)\Big]\sqrt{p} \| \bm{w}\|_2.
    \end{aligned}
\end{equation*}

Using homogeneity, we can assume that $\sqrt{p} \|\bm{w} \|_2 + p^{1 / \theta} \| \bm{w}\|_{\infty} = 1$. Then $\|\bm{w}\|_{2} \leq p^{-1 / 2}$ and $\|\bm{w}\|_{\infty} \leq p^{-1 / \theta}$. Therefore, for $p \geq 2$,
\begin{equation*}
    \begin{aligned}
        \|\bm{w}\|_{p} & \leq \Big(\sum_{i = 1}^n |w_{i}|^{2}\|\bm{w}\|_{\infty}^{p-2}\Big)^{1 / p} \leq (p^{-1-(p-2) / \theta} )^{1 / p}=(p^{-p / \theta} p^{(2-\theta) / \theta })^{1 / p} \\
        & \leq {3^{\frac{2 - \theta}{3 \theta}}}{p^{-1 / \theta}} ={3^{\frac{2 - \theta}{3 \theta}}}{p^{-1 / \theta}} \{\sqrt{p}\|\bm{w}\|_{2}+p^{1 / \theta}\|\bm{w}\|_{\infty}\},
    \end{aligned}
\end{equation*}
where the last inequality follows form the fact that $p^{1 / p} \leq 3^{1 / 3}$ for any $p \geq 2, p \in \mathbb{N}$. Hence
\begin{equation*}
   \begin{aligned}
       \Big\| \sum_{i =  1}^n w_i X_i \Big\|_p & \leq 2 e^{3 + \frac{2 - \theta}{e \theta}} \Gamma^{1 / p} \Big( \frac{p}{\theta} + 1\Big)  \| \bm{w}\|_{\infty} \\
       & \qquad \quad + 2 \bigg[ \log^{1 / \theta} 2 + e^3 \Big( \Gamma^{1 / 2} \Big( \frac{2}{\theta} + 1\Big) + 3^{\frac{2 - \theta}{3 \theta}} p^{-\frac{1}{\theta}}\Gamma^{1 / p} \Big( \frac{p}{\theta} + 1\Big) \Big) \bigg] \sqrt{p} \| \bm{w} \|_2.
   \end{aligned}
\end{equation*}

Following Corollary \ref{lem2}, we have
\begin{equation*}
    \Big\|\sum_{i = 1}^n w_i X_i\Big\|_{\Psi_{\theta, L_n(\theta, p)}} \leq \gamma e D_1 (\theta),
\end{equation*}
where $L_n (\theta, p) = \frac{\gamma^{2 / \theta} D_2 (\theta, p)}{\gamma D_1 (\theta)}$,
$D_1 (\theta) := 2 [ \log^{1 / \theta} 2 + e^3 ( \Gamma^{1 / 2} ( \frac{2}{\theta} + 1) + \sup_{p \geq 2} 3^{\frac{2 - \theta}{3 \theta}} p^{- \frac{1}{\theta}} \Gamma^{1 / p} ( \frac{p}{\theta}$ $+ 1) ) ]\| \bm{w} \|_2 < \infty$,
 and $D_2 (\theta, p) := 2 e^{3} 3^{\frac{2 - \theta}{3 \theta}} p^{- 1 / \theta} \Gamma^{1 / p} \left( \frac{p}{\theta} + 1\right)  \| \bm{w}\|_{\infty}$.

Finally, take $L_n(\theta)  = \inf_{p \geq 1} L_n(\theta, p) > 0.$ Indeed, the positive limit can be argued by (2.2) in \cite{Alzer97}.
Then by the monotonicity property of the GBO norm, it gives
\begin{equation*}
   \Big\|\sum_{i = 1}^n w_i X_i \Big\|_{\Psi_{\theta, L_n(\theta)}} \le \Big\|\sum_{i = 1}^n w_i X_i \Big\|_{\Psi_{\theta, L_n(\theta, p)}} \leq \gamma e D_1 (\theta).
\end{equation*}

\emph{\hl{Case} $\theta > 1$:} In this case $N(t) = t^{\theta}$ is convex with $N^*(t) = \theta^{- \frac{1}{\theta - 1}} \left( 1 - \theta^{-1}\right) t^{\frac{\theta}{\theta - 1}}.$ By Lemmas \ref{lem3} and \ref{lem4}(b), for $p \geq 2$, we have
\begin{equation*}
    \begin{aligned}
      \scalebox{0.92}{$\Big\| \sum_{i =  1}^n w_i Z_i \Big\|_p$} & \scalebox{0.92}{$\leq e \inf \Big\{ t > 0 \, : \, \sum_{i = 1}^n \log \phi_p \Big( \frac{4 w_i}{t} Z_i / 4\Big) \leq p\Big\} + e \inf \Big\{ t > 0 \, : \, \sum_{i = 1}^n p M_{p, Z_i} ( \frac{4 w_i}{t}) \leq p\Big\} $}\\
        & \scalebox{0.92}{$\leq e \inf \Big\{ t > 0 \, : \, \sum_{i = 1}^n p^{-1} N^* \Big( p \Big| \frac{4 w_i}{t}\Big|\Big) \leq 1\Big\} + e \inf \Big\{ t > 0 \, : \, \sum_{i = 1}^n p ( \frac{4 w_i}{t})^2 \leq 1\Big\}$} \\
        & \scalebox{0.92}{$= 4 e \big[  \sqrt{p} \| \bm{w}\|_2 + ( p / \theta)^{1 / \theta} ( 1 - \theta^{-1})^{1 / \beta} \| \bm{w}\|_{\beta}\big]$}
    \end{aligned}
\end{equation*}
with $\beta$ mentioned in the statement. Therefore, for $p \geq 2$, Equation \eqref{eq:symmetrization} implies
\begin{center}
$\| \sum_{i =  1}^n w_i X_i \|_p \leq [ 8 e + 2 (\log 2)^{1 / \theta}] \sqrt{p} \|\bm{w}\|_2 + 8 e ( p / \theta)^{1 / \theta} ( 1 - \theta^{-1})^{1 / \beta} \| \bm{w}\|_{\beta}.$
\end{center}

Then the following result follows by Corollary \ref{lem2},
\begin{center}
    $\left\| \sum_{i = 1}^n w_i X_i \right\|_{\Psi_{\theta, L^{\prime} (\theta)}} \leq \gamma e D_1^{\prime} (\theta) $,
\end{center}
where $L_n (\theta) = \frac{\gamma^{2 / \theta} D_2^{\prime} (\theta)}{\gamma D_1^{\prime} (\theta)}$, $D_1^{\prime}(\theta) = \left[ 8 e + 2 (\log 2)^{1 / \theta}\right]\|\bm{w}\|_2$, and $D_2^{\prime}(\theta) = 8 e \theta^{- 1 / \theta} ( 1 - \theta^{-1})^{1 / \beta} \| \bm{w}\|_{\beta}$.

Note that $w_i X_i = (w_i \|X_i \|_{\psi_{\theta}}) (X_i / \|X_i \|_{\psi_{\theta}})$, we can conclude (a).
\vspace{1ex}

\item \hl{(b)} 
 It is followed from Proposition \ref{prop:GBO norm} and (a).

\vspace{1ex}
\item \hl{(c)} For easy notation, put ${L_n}(\theta) = {L_n}(\theta, \bm{b}_X)$ in the proof. When $\theta  < 2$, by the inequality $a+ b \le 2(a \vee b)$ for $a,b>0$, we have
\begin{center}
    $\pr \left( { | {\sum\limits_{i = 1}^n {{w_i}} {X_i}} | \ge 4eC(\theta ){{\| \bm b \|}_2}\sqrt t } \right) \le 2{e^{ - t}},~~\text{if}~\sqrt t  \ge {L_n}(\theta){t^{1/\theta }}.$
\end{center}

Put $s: = 4eC(\theta ){\left\| \bm b \right\|_2}\sqrt t$, we have
\begin{center}
    $\pr \left( {| {\sum\limits_{i = 1}^n {{w_i}} {X_i}} | \ge s} \right) \le 2\exp \left\{ { - \frac{{{s^2}}}{{16{e^2}{C^2}(\theta )\left\|\bm b \right\|_2^2}}} \right\}, ~\text{ if }~  s \leq 4eC(\theta ){\left\|\bm b \right\|_2} L_n^{\theta /(\theta -  2)}(\theta).$
\end{center}

For $\sqrt t  \le {L_n}(\theta){t^{1/\theta }}$, we obtain $ \pr ( |\sum_{i = 1}^n w_i X_i| \geq 4 e C(\theta) \|\bm{b}_X\|_2 L_n (\theta) t^{ 1 / \theta}) \leq 2 e^{-t}.$ Let $s: = 4eC(\theta ){{\left\| \bm b \right\|}_2}{L_n}(\theta){t^{1/\theta }}$, it gives
\begin{center}
    $\pr \left( {| {\sum\limits_{i = 1}^n {{w_i}} {X_i}} | \ge s} \right) \leq 2\exp \left\{ { - \frac{{{s^\theta }}}{{{{[4eC(\theta ){{\left\|\bm b \right\|}_2}{L_n}(\theta)]}^\theta }}}} \right\},~\text{ if }~s > 4eC(\theta ){\left\|\bm b \right\|_2}L_n^{\theta /(\theta  - 2)}(\theta).$
\end{center}

Similarly, for $\theta  > 2$, it implies
\begin{center}
$ \pr \left({\left| {\sum_{i = 1}^n {{w_i}} {X_i}} \right| \ge s} \right) \le 2e^{ { - \frac{{{s^\theta }}}{{{{[4eC(\theta ){{\left\|\bm  b \right\|}_2}{L_n}(\theta) ]}^\theta }}}} }$ if $s \le 4eC(\theta ){\left\| \bm  b \right\|_2}L_n^{\theta /(2-\theta )}(\theta)$,
\end{center}
and
$\pr \left( {\left| {\sum_{i = 1}^n {{w_i}} {X_i}} \right| \ge s} \right) \le 2e^{ { - \frac{{{s^2}}}{{16{e^2}{C^2}(\theta )\left\|\bm  b \right\|_2^2}}} }$ if $s \ge 4eC(\theta ){\left\|\bm  b \right\|_2}L_n^{\theta /(2-\theta )}(\theta)$.
\end{proof}

\subsection{}

\begin{proof}[\hl{Proof of Corollary} \ref{cor:type1SW}]
Using the definition of $\left\| X \right\|_{{\varphi _\theta }}$, it yields
\begin{align*}
{\rm{E}}{e^{{{({c^{ - 1}}|X|)}^\theta }}} = 1 + \sum\limits_{k = 1}^\infty  {\frac{{{c^{ - k}}{\rm{E}}|X{|^{k\theta }}}}{{k!}}}&  \le 1 + \sum\limits_{k = 1}^\infty  {\frac{{{c^{ - k}}k!\left\| X \right\|_{{\varphi _\theta }}^{k\theta }}}{{k!}}}\\
  &= 1 + \sum\limits_{k = 1}^\infty  {{(\frac{{\left\| X \right\|_{{\varphi _\theta }}^\theta }}{c^{\theta}})^k}}  = 1 + \frac{{\left\| X \right\|_{{\varphi _\theta }}^\theta }}{c^{\theta}}\sum\limits_{k = 0}^\infty  {{(\frac{{\left\| X \right\|_{{\varphi _\theta }}^\theta }}{c^{\theta}})^k}} \\
[\frac{{\left\| X \right\|_{{\varphi _2}}^\theta }}{c^{\theta}} < 1]& = 1 + (\frac{{\left\| X \right\|_{{\varphi _2}}^\theta }}{c^{\theta}})\frac{1}{{1 - \left\| X \right\|_{{\varphi _2}}^\theta /{c^{\theta}}}} \le 2
\end{align*}
if $\frac{{\left\| X \right\|_{{\varphi _2}}^\theta }}{c^{\theta}} \le \frac{1}{2}$ which implies that the minimal $c$ is $2^{1/\theta }\|X\|_{\varphi_{\theta}}$. That is to say we have ${\rm{E}}e^ {|X/[{2^{1/\theta }}\|X\|_{\varphi_{\theta}}]|^{1/\theta}}  \le 2$. Applying \eqref{eq:tail},
we have
\begin{equation}\label{eq: sub-Wcon}
\pr \{ |X| > t \} \le 2{e^{ - {{(t/[{2^{1/\theta }}\|X\|_{{\varphi _\theta }}])}^\theta }}} = 2\exp \{  - \frac{{{t^\theta }}}{{2\left\| X \right\|_{{\varphi _\theta }}^\theta }}\} ~\text{for all } t \ge 0.
\end{equation}
\end{proof}

\subsection{}\label{thm2.proof}

\begin{proof}[\hl{Proof of Theorem} \ref{thm:sub-con}]
Minkowski's inequality for $p\ge 1$ and definition of $\|X\|_{\varphi_{\theta}}$ imply
\[{\left\| {\sum\limits_{i = 1}^n {{X_i}} } \right\|_p} \le \sum\limits_{i = 1}^n {{{\left\| {{X_i}} \right\|}_p}}  \le \sum\limits_{i = 1}^n {{v_i}}  \cdot {2^{1/\theta }}{C_\theta }{\left( {\frac{p}{{\theta {e^{11/12}}}}} \right)^{1/\theta }},\]
where the last inequality by letting ${C_\theta } := \mathop {\max }\limits_{k \ge 1} {\left( {\frac{{2\sqrt {2\pi } }}{\theta }} \right)^{1/k}}{\left( {\frac{k}{\theta }} \right)^{1/(2k)}}$ in Corollary \ref{prop:Psub-W}(b).

From Markov's inequality, it yields
\[\pr\left( {\left| {\sum\limits_{i = 1}^n {{X_i}} } \right| \ge t} \right) \le {t^{ - p}}\left\| {\sum\limits_{i = 1}^n {{X_i}} } \right\|_p^p \le {t^{ - p}}{(\sum\limits_{i = 1}^n {{v_i}} )^p}{2^{p/\theta }}{C_\theta }{\left( {\frac{p}{{\theta {e^{11/12}}}}} \right)^{p/\theta }}.\]

Let ${t^{ - p}}{(\sum\limits_{i = 1}^n {{v_i}} )^p}{2^{p/\theta }}{C_\theta }{\left( {\frac{p}{{\theta {e^{11/12}}}}} \right)^{p/\theta }} = {e^{ - p}}$, it gives
\begin{center}
$t = e(\sum\limits_{i = 1}^n {{v_i}} ){2^{1/\theta }}{C_\theta }{\left( {\frac{p}{{\theta {e^{11/12}}}}} \right)^{1/\theta }}$ and $p = \frac{{\theta {e^{11/12}}{t^\theta }}}{{{[e(\sum\limits_{i = 1}^n {{v_i}} ){2^{1/\theta }}{C_\theta }]^\theta }}}$.
\end{center}

 Therefore, for $p\ge 1$, we have
\begin{equation}\label{sw.con1}
\pr \left( {\left| {\sum\limits_{i = 1}^n {{X_i}} } \right| \ge t} \right) \le \pr \left( {\left| {\sum\limits_{i = 1}^n {{X_i}} } \right| \ge e(\sum\limits_{i = 1}^n {{v_i}} ){C_\theta }{{({2^{ - 1}}\theta {e^{11/12}})}^{ - 1/\theta }}} \right) \le {e^{ - p}} \in (0,{e^{ - 1}}].
\end{equation}

So
\[\pr \left( \left| \sum\limits_{i = 1}^n {{X_i}}  \right| \ge t \right) \le \exp \bigg\{  - \frac{{\theta {e^{11/12}}{t^\theta }}}{{2{[e(\sum\limits_{i = 1}^n {{v_i}} ){C_\theta }]^\theta }}}\bigg\} ,\qquad t \ge e(\sum\limits_{i = 1}^n {{v_i}} ){C_\theta }{{({2^{ - 1}}\theta {e^{11/12}})}^{ - 1/\theta }}.\]

Let $\bar{v}=\frac{1}{n}\sum\limits_{i = 1}^n {{v_i}}$ and ${e^{ - p}} = :\alpha $. Then
\[\pr \left( {\left| {\frac{1}{n}\sum\limits_{i = 1}^n {{X_i}} } \right| \le e\bar v{2^{1/\theta }}{C_\theta }{{\left( {\frac{{\log ({\alpha ^{ - 1}})}}{{\theta {e^{11/12}}}}} \right)}^{1/\theta }}} \right) \ge 1 - \alpha  \in (1 - {e^{ - 1}},1].\]

For $p< 1$, note that moment monotonicity show that $\left[\mathrm{E}\left(|X|^{p}\right)\right]^{1 /p}$ is a non-decreasing function of $p,$, i.e.,
$$
0<p \leq 1 \Rightarrow\left[\mathrm{E}|X|^{p}\right]^{1 / p} \leq\mathrm{E}|X|.
$$

The $c_{r}$ -inequality implies $\left\| {\sum\limits_{i = 1}^n {{X_i}} } \right\|_p^p \le \sum\limits_{i = 1}^n {\left\| {{X_i}} \right\|_p^p}$. Using Markov's inequality again, we have
\[\pr \left( {\left| {\sum\limits_{i = 1}^n {{X_i}} } \right| \ge t} \right) \le {t^{ - p}}\left\| {\sum\limits_{i = 1}^n {{X_i}} } \right\|_p^p \le {t^{ - p}}\sum\limits_{i = 1}^n {\left\| {{X_i}} \right\|_p^p}  \le {t^{ - p}}\sum\limits_{i = 1}^n ({\rm{E}}|{X_i}|)^p .\]

Put $ {t^{ - p}}\sum\limits_{i = 1}^n ({\rm{E}}|{X_i}|)^p = {e^{ - p}}$ and $t = e{( {\sum\limits_{i = 1}^n ({\rm{E}}|{X_i}|)^p })^{1/p}}.$ Then, we obtain
\begin{equation}\label{sw.con2}
    \pr \left( {\left| {\sum\limits_{i = 1}^n {{X_i}} } \right| \ge e{{( {\sum\limits_{i = 1}^n ({\rm{E}}|{X_i}|)^p })^{1/p}}}} \right) \le {e^{ - p}} \in ({e^{ - 1}},1).
\end{equation}

Combine \eqref{sw.con1} and \eqref{sw.con2}, we obtain for all $t \ge 0$,
    \[\pr \left(| \sum\limits_{i = 1}^n {{X_i}} | \ge e{( {\sum\limits_{i = 1}^n ({\rm{E}}|{X_i}|)^t })^{1/t}} + e(\sum\limits_{i = 1}^n {{v_i}} ){2^{1/\theta }}{C_\theta }{( {\frac{t}{{\theta {e^{11/12}}}}} )^{1/\theta }}\right)\le e^{-t}.\]

This completes the proof.
\end{proof}

\subsection{}\label{thm3.proof}

\begin{proof}[\hl{Proof of Theorem} \ref{eq:L2uppercox}]
Note that for $\forall b \in S^{p-1}$, it yields
\begin{align}\label{eq:ZCH}
 &~~~~{{b^{\top}}\hat{\mathcal{Q}}_{n}(\beta^*)b - {b^{\top}}{\rm{E}}(\hat{\mathcal{Q}}_{n}(\beta^*))b} \ge -\left\| b \right\|\mathop {\max }\limits_{k,j} |[\hat{\mathcal{Q}}_{n}(\beta^*) - {\rm{E}}\hat{\mathcal{Q}}_{n}(\beta^*)]_{kj}|\nonumber\\
 &=-\mathop {\max }\limits_{k,j} \left|\frac{1}{n}\sum\limits_{i = 1}^n {\left[ {\frac{{({Y_i} + k)k{e^{ X_i^{\top}\beta^* }}{X_i}X_i^{\top}}}{{{(k + {e^{ X_i^{\top}\beta^* }})^2}}} - {\rm{E}}\left( {\frac{{({Y_i} + k)k{e^{ X_i^{\top}\beta^* }}{X_i}X_i^{\top}}}{{{(k + {e^{ X_i^{\top}\beta^* }})^2}}}} \right)} \right]_{kj}} \right|.
\end{align}

Consider the decomposition
\[\begin{array}{l}
\frac{1}{n}\sum\limits_{i = 1}^n {\left[ {\frac{{({Y_i} + k)k{e^{ X_i^{\top}\beta^* }}{X_{ik}}{X_{ij}}}}{{{(k + {e^{ X_i^{\top}\beta^* }})^2}}} - {\rm{E}}\left( {\frac{{({Y_i} + k)k{e^{ X_i^{\top}\beta^* }}{X_{ik}}{X_{ij}}}}{{{(k + {e^{ X_i^{\top}\beta^* }})^2}}}} \right)} \right]} \\
 = \frac{1}{n}\sum\limits_{i = 1}^n {\left[ {\frac{{{Y_i}k{e^{ X_i^{\top}\beta^* }}{X_{ik}}{X_{ij}}}}{{{(k + {e^{ X_i^{\top}\beta^* }})^2}}} - {\rm{E}}\left( {\frac{{{Y_i}k{e^{ X_i^{\top}\beta^* }}{X_{ik}}{X_{ij}}}}{{{(k + {e^{ X_i^{\top}\beta^* }})^2}}}} \right)} \right]}  + \frac{k}{n}\sum\limits_{i = 1}^n {\left[ {\frac{{k{e^{ X_i^{\top}\beta^* }}{X_{ik}}{X_{ij}}}}{{{(k + {e^{ X_i^{\top}\beta^* }})^2}}} - {\rm{E}}\left( {\frac{{k{e^{ X_i^{\top}\beta^* }}{X_{ik}}{X_{ij}}}}{{{(k + {e^{ X_i^{\top}\beta^* }})^2}}}} \right)} \right]}
\end{array}\]

For the first term, we have under the ${\mathcal{F}_{\max}}$ with $t=C_{\min}/{4}$
\[\begin{array}{l}
     {\pr}\left( { {\left| {\frac{1}{n}\sum\limits_{i = 1}^n {\left[ {\frac{{{Y_i}k{e^{ X_i^{\top}\beta^* }}{X_{ik}}{X_{ij}}}}{{{{(k + {e^{ X_i^{\top}\beta^* }})}^2}}} - {\rm{E}}\left( {\frac{{{Y_i}k{e^{ X_i^{\top}\beta^* }}{X_{ik}}{X_{ij}}}}{{{{(k + {e^{ X_i^{\top}\beta^* }})}^2}}}} \right)} \right]} } \right| \ge t}, {\mathcal{F}_{\max}}} \right)\\
     \le 2\exp \left\{  - \frac{1}{4}\left(\frac{{{n^2}{t^2}}}{{2 \sum\limits_{i = 1}^n {{{({X_{ik}}{X_{ij}})}^2}} {(\left\| {{Y_i}} \right\|_{{\psi _1}}^{} + |\frac{{\exp (X_i^{\top}{\beta ^*})}}{{\log 2}}|)^2}}} \wedge \frac{{nt}}{{\mathop {\max }\limits_{1 \le i \le n} |{X_{ik}}{X_{ij}}|(\left\| {{Y_i}} \right\|_{{\psi _1}}^{} + |\frac{{\exp (X_i^{\top}{\beta ^*})}}{{\log 2}}|)}}\right)\right\} \\
     \le 2\exp \left\{ { - \frac{1}{4}\left(\frac{{n{t^2}}}{{2 M_X^4{{\log }^{4/\theta }}(\frac{{2np}}{\delta })M_{BX}^2}} \wedge \frac{{nt}}{{M_X^2{{\log }^{2/\theta }}(\frac{{2np}}{\delta })M_{BX}^{}}}\right)} \right\}
\end{array}\]
where we use ${{k{e^{ X_i^{\top}\beta^* }}}}{{{(k + {e^{ X_i^{\top}\beta^* }})^{-2}}}} \le 1$ and the second last inequality is from Corollary \ref{remark2}.

For the second term, by Theorem \ref{thm:SumNewOrliczVex} and ${\left\| {{X_{ik}}{X_{ij}}} \right\|_{{\psi _{\theta /2}}}} \le {\left\| {{X_{ik}}} \right\|_{{\psi _\theta }}}{\left\| {{X_{ij}}} \right\|_{{\psi _\theta }}} \le M_X^2$
we have
\[\begin{array}{l}
{\pr}\left( {\left| {\frac{k}{n}\sum\limits_{i = 1}^n {\left[ {\frac{{k{e^{ X_i^{\top}\beta^* }}{X_{ik}}{X_{ij}}}}{{{(k + {e^{ X_i^{\top}\beta^* }})^2}}} - {\rm{E}}\left( {\frac{{k{e^{ X_i^{\top}\beta^* }}{X_{ik}}{X_{ij}}}}{{{(k + {e^{ X_i^{\top}\beta^* }})^2}}}} \right)} \right]} } \right| \ge t},{\mathcal{F}_{\max}} \right)\\
\le 2\exp \left\{ { - \bigg(\frac{{{t^{\theta/2} }}}{{{{[4eC(\theta/2){{\left\| b \right\|}_2}{L_n}(\theta/2 , b)]}^{\theta/2} }}} \wedge \frac{{{t^2}}}{{16{e^2}{C^2}(\theta/2)\left\| b \right\|_2^2}}\bigg)} \right\}
\end{array}\]
where $ b = (k / n) M_X^2(1, \ldots ,1)^{\top}\in\mathbb{R}^n$.

Assume that $ {b^{\top}}{\rm{E}}(\hat{\mathcal{Q}}_{n}(\beta))b\ge C_{\min}$ for all $b \in S^{p-1}$. Under ${\mathcal{F}_1}$ and ${\mathcal{F}_2}$, it shows that by \eqref{eq:ZCH}: ${b^{\top}}{\rm{E}}$\linebreak $(\hat{\mathcal{Q}}_{n}(\beta))b \ge C_{\min} -\frac{ C_{\min}}{2}=\frac{ C_{\min}}{2}$. Then
\begin{align}\label{eq:ZCH1}
 &~~~~ {\pr}\{\lambda _{\min}(\hat{\mathcal{Q}}_{n}(\beta))  \le \frac{ C_{\min}}{2}\} ={\pr}\Big\{{ b^{\top}}{\rm{E}}(\hat{\mathcal{Q}}_{n}(\beta)) b  \le \frac{ C_{\min}}{2},{\forall  b \in S^{p-1} }\Big\}\\
 & \le{\pr}\Big\{{ b^{\top}}{\rm{E}}(\hat{\mathcal{Q}}_{n}(\beta)) b  \le \frac{ C_{\min}}{2},{\forall  b \in S^{p-1} },{\mathcal{F}_{\max}}\Big\} + \pr({\mathcal{F}_{\max}^c})\nonumber\\
 & \le {\pr}\{{\mathcal{F}_1},{\mathcal{F}_{\max}}\} + \pr \{{\mathcal{F}_2},{\mathcal{F}_{\max}}\} + \pr({\mathcal{F}_{R}^c(n)})\nonumber\\
  &\le 2p^2\exp \left\{ { - \frac{1}{4}\bigg(\frac{{n{t^2}}}{{M_X^4{{\log }^{4/\theta }}(\frac{{2np}}{\delta })M_{BX}^2}} \wedge \frac{{nt}}{{M_X^2{{\log }^{2/\theta }}(\frac{{2np}}{\delta })M_{BX}}}\bigg)} \right\}\nonumber\\
    & +2p^2 \exp \left\{ { - \bigg(\frac{{{t^{\theta/2} }}}{{{{[4eC(\theta/2){{\left\| b \right\|}_2}{L_n}(\theta/2 , b)]}^{\theta/2} }}} \wedge \frac{{{t^2}}}{{16{e^2}{C^2}(\theta/2)\left\| b \right\|_2^2}}\bigg)} \right\} + \pr({\mathcal{F}_{\max}^c}).
\end{align}

Then we have by conditioning on  ${\mathcal{F}_1}\cap{\mathcal{F}_2}$
\begin{center}
$\delta_{n}(\beta):=\frac{3}{2}\|[\hat{\mathcal{Q}}_{n}(\beta)]^{-1} \hat{\mathcal{Z}}_{n}(\beta)\|_{2}\le \frac{3}{C_{\min}}\|\hat{\mathcal{Z}}_{n}(\beta)\|_{2}.$
\end{center}

By ${k}/({k + {e^{ X_i^{\top}\beta^* }}}) \le 1$, Corollary \ref{remark2} implies for any $1 \leq k \leq p$,
\begin{equation}\label{eq:2et}
    \begin{aligned}
        \pr \bigg[ \Big|\sqrt{\frac{p}{n}} \sum\limits_{i = 1}^n {\frac{{k({Y_i} - {e^{ X_i^{\top}\beta^* }}){X_{ik}}}}{{k + {e^{ X_i^{\top}\beta^* }}}}} \Big| > & 2{ \Big(\frac{2tp}{{{n}}}\sum\limits_{i = 1}^n {X_{ik}^2\left\| {{Y_i} - {\rm{E}}{Y_i}} \right\|_{{\psi _1}}^2} \Big)^{1/2}} \\
        & + 2t \sqrt{\frac{p}{n}} \mathop {\max }\limits_{1 \le i \le n}{|X_{ik}|{{\left\| {{Y_i} - {\rm{E}}{Y_i}} \right\|}_{{\psi _1}}}} \bigg]  \le 2{e^{ - t}}.
    \end{aligned}
\end{equation}

Let
\begin{center}
${\lambda _{1n}}(t,X) := 2{\Big(\frac{2tp}{{{n}}}\mathop {\max }\limits_{1 \le k \le n}\sum\limits_{i = 1}^n {X_{ik}^2\left\| {{Y_i} - {\rm{E}}{Y_i}} \right\|_{{\psi _1}}^2} \Big)^{1/2}} + 2t \sqrt{\frac{p}{n}} \mathop {\max }\limits_{1 \le i \le n,1 \le k \le p}( {|X_{ik}^{}|{{\left\| {{Y_i} - {\rm{E}}{Y_i}} \right\|}_{{\psi _1}}}} ).$
\end{center}

We bound ${\mathop {\max }\limits_{1 \le i \le n,1 \le k \le p} |{X_{ik}}| \le {M_X}{\log ^{1/\theta }}(\frac{{2np}}{\delta })}$ and $\mathop {\max }\limits_{1 \le k \le n}\frac{{\rm{1}}}{n}\sum\limits_{i = 1}^n {X_{ik}^2} \le {M_X^2}{\log ^{2/\theta }}(\frac{{2np}}{\delta })$ under the event ${\mathcal{F}_{\max}}$. Note that $M_{BX}=M_X + \frac{B}{{\log 2}}$, then (C.1) and (C.2) gives
\begin{equation*}
    \begin{aligned}
        {\lambda _{1n}}(t,X) & \leq 2{\Big(2tpM_{BX}^2\mathop {\max }\limits_{1 \le k \le p} \frac{{\rm{1}}}{n}\sum\limits_{i = 1}^n {X_{ik}^2} \Big)^{1/2}} + 2 t \sqrt{\frac{p}{n}} \mathop {\max }\limits_{1 \le i \le n,1 \le k \le p} |X_{ik}|M_{BX} \\
        & \leq 2M_{BX} M_X (\sqrt{2 t p} + t \sqrt{p / n}) \log^{1 / \theta} (2 n p / \delta) =:{\lambda_{n}}(t).
    \end{aligned}
\end{equation*}

So, $\pr\bigg(|\sqrt{\frac{p}{n}}\sum\limits_{i = 1}^n {\frac{{k({Y_i} - {e^{ X_i^{\top}\beta^* }}){X_{ik}}}}{{k + {e^{ X_i^{\top}\beta^* }}}}} | > {\lambda _{n}}(t) \bigg) \le 2{e^{ - t}},~k=1,2,\ldots,p.$ Thus \eqref{eq:2et} shows
\[\begin{array}{l}
    \pr \{ \sqrt n {\| {{{\widehat {\cal Z}}_n}({\beta ^*})} \|_2} > {\lambda _{1n}}(t)\}  \le \pr \{ \sqrt n {\| {{{\widehat {\cal Z}}_n}({\beta ^*})}\|_2} > {\lambda _{1n}}(t),{\mathcal{F}_{\max}}\}  + \pr ({\mathcal{F}_{\max}^c}) \\
    \le \pr (\bigcup\limits_{k = 1}^p {\{ \| {\frac{1}{{\sqrt n }}\sum\limits_{i = 1}^n {\frac{{k({Y_i} - {e^{ X_i^{\top}\beta^* }}){X_{ik}}}}{{k + {e^{ X_i^{\top}\beta^* }}}}} } \| > \frac{{{\lambda _{1n}}(t)}}{{\sqrt p }}\} } ) + \pr ({\mathcal{F}_{\max}^c}) \le 2p{e^{ - t}} + \pr ({\mathcal{F}_{\max}^c}) =  {\delta }+\varepsilon_n ,
\end{array}\]
where $t := \log (\frac{{2p}}{{{\delta}}})$. Then $\|\hat{\beta}_{n}-\beta^*\|_{2}\le\delta_{n}(\beta^*)\le \frac{3}{C_{\min}}\|\hat{\mathcal{Z}}_{n}(\beta^*)\|_{2}\le \frac{3{\lambda _{1n}}(t)}{C_{\min}{\sqrt n }}$ via Lemma \ref{norm on net general}. 
 Under ${\mathcal{F}_1} \, \cap \, {\mathcal{F}_2} \, \cap \, {\mathcal{F}_{\max}}$, we obtain
\begin{equation*}
    \|\hat{\beta}_{n}-\beta^*\|_{2} \leq \frac{6 M_{BX} M_X}{C_{\min}} \left[ \sqrt{ \frac{2p}{n} \log \Big( \frac{2p}{\delta}\Big)} + \frac{1}{n} \sqrt{p \log \Big( \frac{2p}{\delta}\Big) } \right] \log^{1 / \theta} \Big( \frac{2 n p}{\delta}\Big).
\end{equation*}
Furthermore, under ${\mathcal{F}_1}\cap{\mathcal{F}_2}\cap{\mathcal{F}_{\max}}$, it gives the condition of $n$: \eqref{eq:CC}.
\end{proof}

\subsection{}\label{thm4.proof}

\begin{proof}[\hl{Proof of Theorem} \ref{NA-BY}]
For convenience, the proof is divided into three steps.

\vspace{2ex}

\textit{Step 1}. \hl{Adopting} 
 the lemma
\begin{lemma}[Computing the spectral norm on a net, Lemma 5.4 in \cite{Vershynin18}] \label{net general}
  Let $\mathbf{B}$ be an $p \times p$ matrix,
  and let $\mathcal{N}_{\varepsilon}$ be an $\varepsilon$-net of  $S^{p-1}$
  for some $\varepsilon \in [0,1)$. Then
\begin{center}
  $
\big\|\mathbf{B}\big\|= \max_{||\textbf{x}||_2 =1}\big\|\mathbf{B}\boldsymbol{x}\big\|_2 =\sup_{\boldsymbol{x} \in S^{p-1}} |\langle {\mathbf{B}\boldsymbol{x},\boldsymbol{x}} \rangle|
  \le (1 - 2\varepsilon)^{-1} \sup_{\boldsymbol{x} \in \mathcal{N}_{\varepsilon}} |\langle {\mathbf{B}\boldsymbol{x},\boldsymbol{x}} \rangle|.
  $
\end{center}
\end{lemma}
Then show that $\| \frac{1}{n}\mathbf{A}^{\top}\mathbf{A}-\mathbf{I}_p \| \le 2 \max_{\boldsymbol{x} \in \mathcal{N}_{1/4}} \big| \frac{1}{n} \|\mathbf{A}\boldsymbol{x}\|_2^2 - 1 \big|$. Indeed, note that $\langle \frac{1}{n} \mathbf{A}^{\top}$\linebreak $\mathbf{A}\boldsymbol{x} - \boldsymbol{x}, \boldsymbol{x} \rangle = \langle \frac{1}{n} \mathbf{A}^{\top}\mathbf{A}\boldsymbol{x}, \boldsymbol{x} \rangle - 1 = \frac{1}{n} \|A\boldsymbol{x}\|_2^2 - 1$. By setting $\varepsilon = 1/4$ in Lemma \ref{norm on net general}, we can obtain:
\begin{equation*}
    \big\| \frac{1}{n}\mathbf{A}^{\top}\mathbf{A}-\mathbf{I}_p \big\| \leq ( 1 - 2 \varepsilon ) ^ { - 1 } \sup _ {\boldsymbol{x} \in \mathcal { N } _ { \varepsilon } } | \langle \frac{1}{n} \mathbf{A}^{\top}\mathbf{A}\boldsymbol{x} -\boldsymbol{x}, \boldsymbol{x} \rangle | = 2 \max _ {  \boldsymbol{x}  \in \mathcal { N } _ { 1 / 4 } } \Big| \frac { 1 } { n } \|  \mathbf{A} \boldsymbol{x}  \| _ { 2 } ^ { 2 } - 1 \Big|.
\end{equation*}

\textit{\hl{Step 2}}. Let $Z_i:=|\langle {\boldsymbol{A}_i,\boldsymbol{x}} \rangle|$ fix any $\boldsymbol{x} \in S^{n-1}$. Observe that
$\|\mathbf{A}\boldsymbol{x}\|_2^2 = \sum_{i=1}^n |\langle {\boldsymbol{A}_i,\boldsymbol{x}} \rangle|^2 =\sum_{i=1}^n Z_i^2$. The fact that $\{ {Z_i}\} _{i = 1}^n$ are $\operatorname{subW}(\theta)$ with ${\rm E} Z_i^2 = 1, \max_{1 \le i \le n}\|Z_i\|_{\psi _\theta}= K$. Then by Corollary \ref{lem:2theta}, $Z_i^2 $ are independent $\operatorname{subW}(\theta/2)$ r.v.s with $\max_{1 \le i \le n}\|Z_i^2\|_{\psi _{\theta/2}}= K^2$. The norm triangle inequality (Lemma A.3 in \cite{Gotze2019}) gives
\begin{equation}\label{eq:bounds}
    \mathop {\max }\limits_{1 \le i \le n}\|Z_i^2-1\|_{\psi _{\theta/2}}\le K_{\theta/2}[1+([({e\theta/2} )^{\theta/2}] \log 2)^{-\theta/2}]K.
\end{equation}
where $K_{\alpha}:=2^{1 / \alpha}$ if $\alpha \in(0,1)$ and $K_{\alpha}=1$ if $\alpha \geq 1$.

Denote $\bm{b}_X := \frac{1}{n} (\|Z_1^2 - 1\|_{\psi _{\theta/2}}, \ldots, \|Z_n^2 - 1\|_{\psi _{\theta/2}})^{\top}$ in Theorem \ref{thm:SumNewOrliczVex}. With \eqref{eq:bounds}, we have
\begin{center}
${\left\| {\bm{b}_X} \right\|_2} ={n^{-1}}\sqrt{\sum_{i = 1}^n {{\| {{Z_i^2}-1} \|^2_{{\psi _{\theta/2} }}}}}\le \frac{K_{\theta/2}[1+([({e\theta/2} )^{\theta/2}] \log 2)^{-\theta/2}]K}{{\sqrt{n}}}$
\end{center}
and $\| \bm b\|_{\infty}\le \frac{{K_{\theta/2}[1+([({e\theta/2} )^{\theta/2}] \log 2)^{-\theta/2}]K}}{{{n}}}$.

For $\beta : = \frac{\theta }{{\theta  - 1}}>1$, we obtain
\begin{center}
${\left\| {\bm{b}_X} \right\|_{\beta}} ={n^{-1}}
\{\sum_{i = 1}^n {{\| {{Z_i^2}-1} \|^{\beta}_{{\psi _{\theta/2} }}}}\}^{1/\beta}\le {n^{{\beta ^{ - 1}} - 1}}[{K_{\theta/2}[1+([({e\theta/2} )^{\theta/2}] \log 2)^{-\theta/2}]K}]={n^{{-\theta ^{ - 1}} }}{K_{\theta/2}[1+([({e\theta/2} )^{\theta/2}] \log 2)^{-\theta/2}]K}$.
\end{center}

Write ${L_n}(\theta {\rm{/2}},{\bm{b}_X})$ as the constant defined in Theorem \ref{thm:SumNewOrliczVex}(a). Then,
\begin{align*}
    \| \bm{b}_X\|_2 L_n ({\theta}/{2},\bm{b}_X)& =\gamma^{4 / \theta} \left\{ \begin{array}{ll}
           A(\theta / 2) \| \bm b\|_{\infty},&\theta \leq 2\\
           B(\theta / 2) \| \bm b\|_{\beta},&\theta>2.
        \end{array} \right. \\
        &\leq {K_{\theta/2}[1+([({e\theta/2} )^{\theta/2}] \log 2)^{-\theta/2}]K}\gamma^{4 / \theta} \left\{ \begin{array}{ll}
            A(\theta / 2) /n, & \theta \leq 2 \\
            B(\theta / 2)/{n}^{1/\theta }, &   \theta>2 .
        \end{array} \right.
\end{align*}

Hence
\begin{align*}
&~~~~2eC(\theta/2 )\{{\left\|\bm{b}_X \right\|_2} \sqrt t + {\left\|\bm b \right\|_2}{L_n}(\theta/2 ,{\bm{b}_X}){t^{2/\theta }}\} \\
        &  \leq 2eKC(\theta/2 ){K_{\theta/2}[1+([({e\theta/2} )^{\theta/2}] \log 2)^{-\theta/2}}) \left[ \sqrt{\frac{t}{n}} + \left\{ \begin{array}{ll}
            A(\theta / 2) (\gamma^2 t)^{2 / \theta} / {n},  &  \theta \leq 2 \\
            B(\theta / 2) (\gamma^2 t)^{2 / \theta} /{n}^{1/\theta }, & \theta>2
        \end{array} \right. \right]\\
        &=: H(t, n ; \theta).
\end{align*}

Therefore, $\pr ( \frac{1}{n} | \sum_{i = 1}^n (Z_i^2  - 1) | \geq H (t, n ; \theta ) ) \leq 2 e^{-t}$. Let $t = c p + s^2$ for constant $c$, then
\begin{equation*}
    {\pr}\bigg\{ \Big| \frac {1} {n} \|  \mathbf{A} \boldsymbol{x}  \| _ { 2 } ^ { 2 } - 1 \Big| \geq H(cp + s^2, n ; \theta)\bigg\} \leq 2e^{-(c p+s^2)}.
\end{equation*}

\textit{\hl{Step 3}}. Consider the follow lemma for covering numbers in \cite{Vershynin18}.

\begin{lemma}[Covering numbers of the sphere]\label{net cardinality}
 For the unit Euclidean sphere $S^{n-1}$, the covering number $\mathcal{N}(S^{n-1},\varepsilon)$ satisfies $\mathcal{N}(S^{n-1},\varepsilon) \le ( 1 + \frac{2}{\varepsilon} )^n $ for every $\varepsilon >0$.
\end{lemma}
Then, we show the concentration for $\| \frac{1}{n}\mathbf{A}^{\top}\mathbf{A}-\mathbf{I}_p \| $, and \eqref{eq:P2} follows by the definition of largest and least eigenvalues. The conclusion is drawn by Step 1 and 2:
\begin{align*}
&~~~~\pr \bigg\{ \Big\| \frac{1}{n}\mathbf{A}^{\top}\mathbf{A}-\mathbf{I}_p \Big\| \geq H(cp + s^2, n; \theta) \bigg\}  \leq \pr \bigg\{ 2 \max_{  \bm{x}  \in \mathcal {N}_{1/4} } \Big| \frac{1}{n} \| \mathbf{A} \boldsymbol{x}\|_{ 2 }^{ 2 } - 1\Big|\geq H(cp + s^2, n; \theta) \bigg\} \\
        & \leq \mathcal{N} (S^{n - 1}, 1 / 4) \pr \bigg\{  \Big| \frac{1}{n} \| \mathbf{A} \boldsymbol{x}\|_{ 2 }^{ 2 } - 1\Big| \geq H(cp + s^2, n; \theta) / 2 \bigg\}\leq 2 \cdot 9^n e^{-(cp + s^2)},
\end{align*}
where the last inequality follows by Lemma~\ref{net cardinality} with $\varepsilon = 1/4$. When the $c \geq {n \log 9/p}$, then $2 \cdot 9^n e^{-(c p+s^2)} \leq 2 e^{-s^2}$, and the \eqref{eq:P2} is proved.

Moreover, note that
\begin{equation*}
    \max_{||\boldsymbol{x}||_2 =1} \Big| \big\|\frac{1}{\sqrt n}\mathbf{A}\boldsymbol{x}\big\|_2^2-1 \Big| = \max_{||\boldsymbol{x}||_2 =1} \big\|(\frac{1}{n}\mathbf{A}^{\top}\mathbf{A}-\mathbf{I}_p)\boldsymbol{x}\big\|_2^2 = \big\| \frac{1}{n}\mathbf{A}^{\top}\mathbf{A}-\mathbf{I}_p \big\|^2\leq H^2(c p + s^2, n ; \theta).
\end{equation*}
implies that
\begin{center}
$\sqrt{1 - H^2(c p + s^2, n ; \theta)} \le \frac{1}{\sqrt{n}} \lambda_{max}(\mathbf{A}) \le \sqrt{1+H^2(c p + s^2, n ; \theta)}$.
\end{center}

Similarly, for the minimal eigenvalue, we have
\begin{equation*}
    \min_{||\boldsymbol{x}||_2 =1} \Big| \big\|\frac{1}{\sqrt n}\mathbf{A}\boldsymbol{x}\big\|_2^2-1 \Big| = \min_{||\boldsymbol{x}||_2 =1} \big\|(\frac{1}{n}\mathbf{A}^{\top}\mathbf{A}-\mathbf{I}_p)\boldsymbol{x}\big\|_2^2 = \big\| \frac{1}{n}\mathbf{A}^{\top}\mathbf{A}-\mathbf{I}_p \big\|^2\leq H^2(c p + s^2, n ; \theta).
\end{equation*}

This implies $\sqrt{1 - H^2(c p + s^2, n ; \theta)} \le \frac{1}{\sqrt{n}} \lambda_{min}(\mathbf{A}) \le \sqrt{1+H^2(c p + s^2, n ; \theta)}$. So we obtain that the two events satisfy
\begin{equation*}
    \begin{aligned}
        & \Big\{\big\| \frac{1}{n}\mathbf{A}^{\top}\mathbf{A}-\mathbf{I}_p \big\|^2\leq  H^2(c p + s^2, n ; \theta) \Big\} \\
        & \qquad \subset \Big\{\sqrt{1- H^2(c p + s^2, n ; \theta)} \le \frac{1}{\sqrt{n}}\lambda_{min}(\mathbf{A}) \le \frac{1}{\sqrt{n}}\lambda_{max}(\mathbf{A}) \le \sqrt{1 + H^2(c p + s^2, n ; \theta)}\Big\}
    \end{aligned}
\end{equation*}

Then we obtain the second conclusion in this theorem.
\end{proof}

\subsection{}
\begin{proof}[\hl{Proof of Theorem} \ref{eq:cantoniG}]
By independence and \eqref{eq:CANTON},
\vspace{-6pt}

\begin{align}\label{eq;mgfbound}
 &\mathrm{E}{e^{ \pm n\alpha_n {{\hat Z}_{\alpha_n} }(\theta )}}=  \prod\limits_{i = 1}^n {\rm{E}}\exp \{  \pm \varphi^c \left[ {\alpha_n \left( {{X_i} - {\theta }} \right)} \right]\}  \le \prod\limits_{i = 1}^n {\rm{E}}[1 \pm \alpha_n \left( {{X_i} - {\theta }} \right) + c( {\alpha_n ( {{X_i} - {\theta }} )} )]\nonumber\\
 &\le \prod\limits_{i = 1}^n \exp \{  \pm \alpha_n {\rm{E}}\left( {{X_i} - {\theta }} \right) + {\rm{E}}[c( {\alpha_n ( {{X_i} - {\theta }} )} )] \}= \exp \left\{  \pm \alpha_n \sum\limits_{i = 1}^n {\rm{E}}\left( {{X_i} - {\theta }} \right)  + \sum\limits_{i = 1}^n {{\rm{E}}[c( {\alpha_n ( {{X_i} - {\theta }} )} )] }  \right\}.
\end{align}

For convenience, let
	\begin{align}\label{eq:Bplus}
		B_n^ + (\theta ) = \mu_n  - \theta  + \frac{1}{{n\alpha_n }}\sum\limits_{i = 1}^n {{\rm{E}}[c\left( {\alpha_n ({X_i} - \theta) } \right)]} {\rm{ + }}\frac{{\log ({\delta ^{ - 1}})}}{n\alpha_n }
	\end{align}
	and
	$B_n^ - (\theta ) = \mu_n  - \theta  - \frac{1}{{n\alpha_n }}\sum\limits_{i = 1}^n {{\rm{E}}[c\left( {\alpha_n ({X_i} - \theta )} \right)]}  - \frac{{\log ({\delta ^{ - 1}})}}{n\alpha_n }$.	
	Therefore, Equation \eqref{eq;mgfbound} and the Markov's  inequality show
	\vspace{-6pt}
	\[\pr ({{\hat Z}_{\alpha_n} }(\theta ) \ge B_n^ + (\theta )) = \pr ({e^{n\alpha_n {{\hat Z}_{\alpha_n} }(\theta )}} \ge {e^{n\alpha_n B_n^ + (\theta )}}) \le \frac{{{\rm{E}}{e^{n\alpha_n {{\hat Z}_{\alpha_n} }(\theta )}}}}{{{e^{n\alpha_n B_n^ + (\theta )}}}} \le \frac{{{e^{n\alpha_n B_n^ + (\theta ) - \log ({\delta ^{ - 1}})}}}}{{{e^{n\alpha_n B_n^ + (\theta )}}}} = \delta \]
	and
	$\pr (  {{\hat Z}_{\alpha_n} }(\theta ) \le  B_n^ - (\theta )) = \pr({e^{ - n\alpha_n {{\hat Z}_{\alpha_n} }(\theta )}} \ge {e^{ - n\alpha_n B_n^ - (\theta )}}) \le \frac{{{\rm{E}}{e^{ - n\alpha_n {{\hat Z}_{\alpha_n} }(\theta )}}}}{{{e^{ - n\alpha_n B_n^ - (\theta )}}}} \le \frac{{{e^{ - n\alpha_n B_n^ - (\theta ) - \log ({\delta ^{ - 1}})}}}}{{{e^{ - n\alpha_n B_n^ - (\theta )}}}} = \delta $.	
These two inequality  yield $\pr \big(B_n^ - (\theta ) < {{\hat Z}_{\alpha_n} }(\theta ) \big) = 1- \pr \big(  {{\hat Z}_{\alpha_n} }(\theta ) \le  B_n^ - (\theta )\big) - \pr \big({{\hat Z}_{\alpha_n} }(\theta ) \ge B_n^ + (\theta ) \big) \ge 1 - 2\delta .$\\

	The $\frac{{\partial {{\hat Z}_{\alpha_n} }(\theta )}}{{\partial \theta }} = -\frac{{ 1}}{n}\sum_{i = 1}^n {\dot \varphi^c } \left[ {\alpha_n \left( {{X_i} - {\theta }} \right)} \right] < 0$ implies the map $\theta \mapsto {{\hat Z}_{\alpha_n} }(\theta )$ is non-increasing.	
	If $\theta = \mu_n$, we have $B_n^ + (\mu_n)>0$ from \eqref{eq:Bplus}. As $n$ is sufficient large and $\alpha_n \to 0$, in $B_n^ + (\theta )$, from (C.1.2) the term $\frac{1}{{n\alpha_n }}\sum_{i = 1}^n {{\rm{E}}[c\left( {\alpha_n ({X_i} - \theta )} \right)]} \le \frac{f(\alpha_n)}{{\alpha_n }}\frac{1}{{n}}\sum_{i = 1}^n {{\rm{E}}[c({X_i} - \theta )]}=\frac{f(\alpha_n)}{{\alpha_n }}O(1)$ converges to 0 by (C.1.3). Then, there must be a constant $d_n(c)>0$ such that $B_n^ + (\mu_n+d_n(c) )<0$. So under \eqref{eq:sample condition}, it implies that $B_n^ + (\theta )=0$ has a solution and denote the smallest solution $\theta_{+} \in (\mu_n,\mu_n+d_n(c))$.	
	Similarly, for $B_n^ - (\theta )$, we have $B_n^ - (\mu_n)<0$. The condition \eqref{eq:sample condition} implies $B_n^ - (\mu_n-d_n(c) )>0$, then $B_n^ - (\theta )=0$ has a solution and denote the largest solution  $\theta_{-} \in (\mu_n-d_n(c),\mu_n)$.	
	Note that ${{\hat Z}_{\alpha_n} }(\theta )$ is a continuous and non-increasing  function, the estimating equation ${{\hat Z}_{\alpha_n} }(\theta ) =0$ has a solution $\hat{\theta}_{\alpha_n} \in [\theta_{-},\theta_{+}]$ such that
	$\theta_{-} \leq \hat{\theta}_{\alpha_n} \leq \theta_{+} $ with a probability at least $1-2\delta$.	
	Recall that
	\begin{equation}
		B_n^ + ({\theta _{\rm{ + }}}) = \mu_n  - {\theta _{\rm{ + }}} + \frac{1}{{n\alpha_n }}\sum\limits_{i = 1}^n {{\rm{E}}[c\left( {\alpha_n ({X_i} - {\theta _{\rm{ + }}})} \right)]}  +\frac{{\log ({\delta ^{ - 1}})}}{n\alpha_n } = 0.\label{eq: negativeB}
	\end{equation}
has the smallest solution $\theta_{+} \in (\mu_n,\mu_n+d_n(c))$ under the condition \eqref{eq:sample condition}. We have
	\begin{equation}
		\begin{aligned}
			\mu_n  - \hat \theta_{\alpha_n}  \ge \mu_n  - {\theta _{\rm{ + }}}&=  -\frac{1}{{n\alpha_n }}\sum\limits_{i = 1}^n {{\rm{E}}[c\left( {\alpha_n {X_i} - \alpha_n {\theta _{\rm{ + }}}} \right)]}  - \frac{{\log ({\delta ^{ - 1}})}}{n\alpha_n }\\
			&=  - \frac{1}{{n\alpha_n }}\sum\limits_{i = 1}^n {\rm{E}}[c\left( {\alpha_n ({X_i} -  \mu_n) } +  {\alpha_n (\mu_n  - {\theta _{\rm{ + }}})} \right)]  - \frac{{\log ({\delta ^{ - 1}})}}{n\alpha_n }\\
[\text{By (C.1.1)}]~	&\ge   - \frac{{{c_2}}}{{n\alpha_n }}\sum\limits_{i = 1}^n {{\rm{E}}[c\left( {\alpha_n {X_i} - \alpha_n \mu_n } \right)]}  - \frac{{{c_2}}}{\alpha_n }\cdot c\left( {\alpha_n (\mu_n  - {\theta _{\rm{ + }}})} \right) - \frac{{\log ({\delta ^{ - 1}})}}{n\alpha_n }
		\end{aligned}
		\end{equation}
		which implies
		\begin{align}\label{eq:solutions}
			\mu_n  - {\theta _{\rm{ + }}} +\frac{{{c_2}}}{\alpha_n }\cdot c\left( {\alpha_n (\mu_n  - {\theta _{\rm{ + }}})} \right) \ge  - \left( {\frac{{{c_2}}}{{n\alpha_n }}\sum\limits_{i = 1}^n {{\rm{E}}[c\left( {\alpha_n ({X_i} - \mu_n )} \right)] + \frac{{\log ({\delta ^{ - 1}})}}{n\alpha_n }} } \right).
		\end{align}

Put $\frac{{{c_2}}}{{n\alpha_n }}\sum_{i = 1}^n {{\rm{E}}[c\left( {\alpha_n ({X_i} - \mu_n )} \right)] = \frac{{\log ({\delta ^{ - 1}})}}{n\alpha_n }}$, i.e., $\sum_{i = 1}^nc_2 {{\rm{E}}[c\left( {\alpha_n ({X_i} - \mu_n )} \right)]}  = \log ({\delta ^{ - 1}}).$ The scaling assumption $c(tx) \le f(t)c(x)$ gives $$f(\alpha_n )c_2\sum_{i = 1}^n {{\rm{E}}[c\left({{X_i} - \mu_n } \right)]}  \ge c_2\sum_{i = 1}^n {{\rm{E}}[c\left( {\alpha_n ({X_i} - \mu_n )} \right)]}  = \log ({\delta ^{ - 1}})$$
and thus $\alpha_n  \ge {f^{ - 1}}\left( {\frac{{\log ({\delta ^{ - 1}})}}{{c_2\sum\nolimits_{i = 1}^n {{\rm{E}}[c {({X_i} - \mu_n )} ]} }}} \right)$.	
	Let ${g_{\alpha_n} }(t) = t + \frac{{{c_2}}}{\alpha_n }c\left( {\alpha_n t} \right)$. Moreover, \mbox{Equation \eqref{eq:solutions}} and the value $\alpha_n$ yields
	\vspace{-6pt}
	\[{g_{\alpha_n} }(\mu_n  - {\theta _{\rm{ + }}} ) =\mu_n  - {\theta _{\rm{ + }}} + \frac{{{c_2}}}{\alpha_n }c\left( {\alpha_n (\mu_n  - {\theta _{\rm{ + }}})} \right) \ge  -{\frac{{2\log ({\delta ^{ - 1}})}}{n\alpha_n }}.\]
	
	Solve the above inequality in terms of $\mu_n  - {\theta _{\rm{ + }}}$, we obtain
\begin{center}
$\mu_n  - \hat \theta_{\alpha_n}  \ge \mu_n  - {\theta _{\rm{ + }}} \ge g_{\alpha_n} ^{ - 1}\left\{-{\frac{{2\log ({\delta ^{ - 1}})}}{n\alpha_n }}\right\}.$
\end{center}

	Similarly, for ${\theta _{\rm{- }}}$, one has
	$\mu_n  - \hat \theta_{\alpha_n}  \le \mu_n  - {\theta _ - } \le  - g_{\alpha_n} ^{ - 1}\left\{-{\frac{{2\log ({\delta ^{ - 1}})}}{n\alpha_n }}\right\}$. Then we obtain that \eqref{eq:con} holds with probability at least $1-2\delta$.
\end{proof}


\end{document}